\documentclass[12pt]{article}
\usepackage[english]{babel}
\usepackage{amsmath}
\usepackage{amssymb}
\usepackage{amsfonts}
\usepackage{amsthm}
\usepackage{graphics,graphicx}
\usepackage{color}
\theoremstyle{plain}

\newtheorem{theorem}{Theorem}
\newtheorem*{proposition*}{Proposition}

\newtheorem*{corollary*}{Corollary}
\newtheorem*{Corollary*}{Important corollary}
\newtheorem{remark}{Remark}

\theoremstyle{definition}

\theoremstyle{remark}

\newtheorem*{remark*}{Remark}

\newtheorem*{example*}{Example}

\newtheorem*{examples*}{Examples}

\numberwithin{equation}{section}
\numberwithin{lemma}{section}
\numberwithin{theorem}{section}
\numberwithin{hypothesis}{section}
\numberwithin{definition}{section}
\numberwithin{example}{section}
\numberwithin{corollary}{section}
\numberwithin{remark}{section}
\begin{document}

\begin{center}
{\large\bf Factorization of a class of matrix-functions \\ with stable partial indices}

\vspace{3mm}
{\bf Gennady MISHURIS$^1$, Sergei ROGOSIN$^{1, 2}$}

\vspace{3mm}
{\footnotesize\it $^1$Aberystwyth University, Penglais, SY23 3BZ Aberystwyth, UK;

\vspace{1mm}
e-mails: ggm@aber.ac.uk; ser14@aber.ac.uk}

\vspace{3mm}
{\footnotesize\it $^2$Belarusian State University, Nezavisimosti Ave., 4, 220030 Minsk, Belarus;

\vspace{1mm}
e-mail: rogosin@bsu.by}

\end{center}

{\footnotesize {\bf Abstract.}
{A new effective method for factorization of a class of nonrational  $n\times n$ matrix-functions with
\emph{stable partial indices} is proposed. The method is a generalization of the one recently proposed by the authors which was
valid for the canonical factorization only. The class of considered matrices is motivated by problems originated
from applications. The properties and details of the asymptotic procedure are illustrated
by examples. The efficiency of the procedure is highlighted by numerical results. }



\section{Introduction}

We {consider} here the problem of factorization of continuous
matrix functions of the real variable. {This} means the representation of a given invertible square matrix
$G \in\left( {\mathcal C}({\mathbb R})\right)^{n\times n}$ in the following form
\begin{equation}
\label{fact}
G(x) = G^{-}(x) \Lambda(x)  G^{+}(x),
\end{equation}
where continuous invertible matrices $G^{-}(x)$, $G^{+}(x)$ possess an analytic continuation in the lower
$\Pi^{-} = \{z = x + i y: {\mathrm{Im}}\, z < 0\}$ and upper $\Pi^{+} = \{z = x + i y: {\mathrm{Im}}\, z > 0\}$
half-planes, respectively, and
\begin{equation}
\label{fact1}
\Lambda(x) = {\mathrm{diag}}\, \left(\left(\frac{x - i}{x + i}\right)^{{\ae}_1}, \ldots, \left(\frac{x - i}{x + i}\right)^{{\ae}_n}\right),\; {\ae}_1, \ldots, {\ae}_n\in {\mathbb Z}.
\end{equation}
The representation (\ref{fact}) is called {\it right (continuous or standard) factorization} and can be considered for any oriented curve $\Gamma$ of {a} certain classes which divides {the} complex plane into two domain $D^{-}$ and $D^{+}$ with {a} changing of diagonal entries in $\Lambda(x)$ for $\left(\frac{x - t^{+}}{x - t^{-}}\right)^{{\ae}_j}$, $t^{\mp}\in D^{\mp}$, or for $x^{{\ae}_j}$ {(if $0\in D^{+}$)}. A similar representation
$$
G(x) = G^{+}(x) \Lambda(x)  G^{-}(x)
$$
is called {\it left (continuous or standard) factorization}. If {the} right- (left-) factorization exists{,} then the integer numbers ${\ae}_1, \ldots, {\ae}_n$, called {\it partial indices}, are determined uniquely up to the order. {In particular, there exist constant transformations of factors, such that ${\ae}_1 \geq \ldots \geq {\ae}_n$.} The factors $G^{-}$, $G^{+}$ are not unique. Relations of pairs of factors are described, e.g. in \cite{LitSpi87}.  The right- (left-) factorization is called {\it canonical factorization} if all partial indices are equal to 0, i.e. ${\ae}_1 = \ldots = {\ae}_n = 0$.

Factorization of matrix functions was first studied in relation to the vector-matrix Riemann (or Riemann-Hilbert) boundary value problem (see \cite{Gak77}). The later was formulated by Riemann in his work on construction of complex differential equations with algebraic coefficients having a prescribed monodromy group (see, e.g., \cite{EhrSpi01}). By {using the method of the Cauchy type integral,} the  vector-matrix Riemann boundary value problem was reduced in \cite{Mus68}, \cite{Vek67}  to a system of the Fredholm integral equations.  A part of the theory of the factorization problem is based on the study of such systems (see also \cite{Gak52}), though this approach does not answer, in particular, the questions of when it is possible to get factorization, how to construct factors and how to determine partial indices.

Among other sources of interest to the factorization problems one can point out the vector-valued Wiener-Hopf equations on a half-line (see \cite{GohKre58}, \cite{Spe85}, \cite{Nob88}) and their discrete analog{ie}s, namely the block Toeplitz equations (see, e.g., \cite{BoeSil90},  \cite{GohFel74}).  The developed technique found several applications in diffraction theory, fracture mechanics, geophysics, financial mathematics etc. (see a brief description given, e.g. in \cite{VeiAbr07} and references therein).

Theoretical background for the study of the matrix factorization and its numerous generalizations is presented in \cite{BoeSpi13}, \cite{EhrSpi01}, \cite{LitSpi87}, \cite{Vek67} (see also \cite{LawAbr07}). The theory of the factorization is more or less complete (see \cite{BoeSpi13}), but the above mentioned constructive questions about existence, factors and partial indices (which are very important for practical applications) have been answered only in a number of special cases.  Among them one can mention rational matrix functions (see, e.g., \cite{Gak52}), functional commutative matrix functions (those satisfying $G(t) G(s) = G(s) G(t), \forall t, s\in \Gamma$, see \cite{Che56a}), upper- (lower-) triangular matrices with factorizable diagonal elements (see \cite{Che56}, \cite{FelGohKru94}), certain classes of meromorphic matrix functions (see \cite{Adu92}, \cite{AmiKam07}, \cite{Kiy12}), special cases of $2\times 2$  Daniele-Khrapkov matrix functions (with a small degree of deviator polynomial) (see \cite{Dan78}, \cite{EhrSpe02}, \cite{Khr71a}, \cite{Khr71b}), special cases of $2\times 2$ matrix functions with three rationally independent entries (see \cite{Akt92}, \cite{EhrSpe02}{, \cite{ProSpe90}} and references therein), special cases of $n\times n$ generalization of  the Daniele-Khrapkov matrix functions (see, e.g., \cite{CaSaMa01}, \cite{Jon84a}, \cite{Raw85}, \cite{VeiAbr07}), special classes of matrices possessing certain symmetry property (see \cite{Vor11}  and references therein). Some approximate and asymptotic methods for matrix factorization have been developed too (see \cite{Cri01}, \cite{JanLagEph99}, \cite{Koi54}). For more detailed description of constructive approaches to solve factorization problem we refer to the recent survey by the authors \cite{MishRog15a} and references therein.

In \cite{MishRog14a} we propose a new asymptotic method of construction of factors for a special class of  $n\times n$ nonrational matrix functions
in the case canonical factorization. The essential property of the considered matrices is that they become close (after suitable transformation) to a unit matrix ({a} similar assumption is used in \cite{Akt92}, \cite{Cam_dSan00}). The idea to use such representation in factorization is related to that in general operator theory and has been exploited since the seminal work by Gohberg and Krein \cite{GohKre58} (see also \cite{ClaGoh87}, \cite{LitSpi87}).  By reduction of the factorization problem to a special matrix boundary value problem we succeed in \cite{MishRog14a} to get at each step of approximation a solution retaining the main properties of the given matrix. It helps to find conditions under which obtained asymptotic series converge.

In the present paper we generalize the approach of \cite{MishRog14a} and apply a new more general procedure in the case stable partial indices.
Stability of partial indices is of special interest for factorization theory (see \cite{Boj58}, \cite{GohKre58b}, and recent paper \cite{Kis15}) as it justifies the use of approximate procedures for factorization. In the considered case the stability of partial indices allows us to reduce auxiliary relations for components of our asymptotic formula to $n^2$ independent scalar Riemann boundary value problems. The crucial point is that some of these problems have negative index equal to $-1$. Thus, their vanishing at infinity solutions exist if and only if the given matrix function satisfies certain solvability conditions. Instead, we consider all above scalar Riemann boundary value problems in the class of functions bounded at infinity. It leads to corrections of the  algorithm of \cite{MishRog14a}. Finally we get both factors of factorization problem in form of asymptotic series. The conditions of their convergence are found too.  

To the best of authors' knowledge, {our} class does not coincide with any of the above mentioned classes. This class contains matrix functions {which appears in} the study of certain problems in fracture mechanics related to perturbation of the crack propagation (\cite{AMS14}, \cite{MishMovMov10}, \cite{PiMiMo07} -- \cite{PiMiMo12}). {Another motivation is the use of such matrix functions in the study the inverse scattering problem (see \cite{Cam_dSan00}).

The paper is organized as follows. In Sec. 2 we introduce necessary notation and formulate the problem.
A constructive algorithm is presented in Sec. 3.   The method is illustrated by an example given in Sec. 4. We conclude our study by showing the quality of the factorization approximation by restricting ourselves only
to the first asymptotic term and discuss the role of the chosen small parameter.

\section{A class of matrices. Problem formulation}
\label{class}

Let us introduce the following class of invertible continuous $n\times n, n \geq 2$, matrix-functions ${\mathcal S} K_n$
depending on a real parameter $\varphi\in {\mathbb R}$,
satisfying the following conditions:

\noindent (1) $G_{\varphi}\in \left({\mathcal C}({\mathbb R})\right)^{n\times n}$ belongs to ${\mathcal S} K_n$ if it can be represented in the form
\begin{equation}
\label{cond1}
G_{\varphi} = R_{\varphi} F R_{\varphi}^{-1},
\end{equation}
where bounded locally H\"older-continuous on ${\mathbb R}$ (in general non-rational) invertible matrix $R_{\varphi}$ is such that

\noindent (2)
\begin{equation}
\label{cond2}
R_{0} = R_{\varphi}\bigl|_{\varphi=0} = I,
\end{equation}

\noindent (3) matrix function $F$ does not depend on parameter $\varphi$, has H\"older-continuous entries $f_{kl}\in H_{\mu}(\overline{\mathbb R})$, i.e. for all $k,l = 1, \ldots, n,$
\begin{equation}
\label{cond3}
\left|f_{kl}(x_1) - f_{kl}(x_2)\right| \leq C \left|\frac{1}{x_1 + i} - \frac{1}{x_2 + i}\right|^{\mu},\; \forall x_1, x_2\in \overline{\mathbb R}, \;
0 < \mu < 1,
\end{equation}
and there exist
$$
\lim\limits_{x\rightarrow\infty} f_{kl}(x) = f_{kl}(\infty) = f_{kl}(-\infty) = f_{kl}(+\infty),
$$
satisfies the following asymptotic estimate at infinity

\noindent (4)
\begin{equation}
\label{cond4}
F(x) \rightarrow I,\;\;\; |x| \rightarrow \infty,
\end{equation}

\noindent (5) $F$ admits a right  factorization, i.e.
\begin{equation}
\label{cond5}
F(x) = F^{-}(x) \Lambda(x) F^{+}(x),
\end{equation}
where H\"older-continuous on $\overline{\mathbb R}$ matrix-functions $F^{-}(x), F^{+}(x)$ possess an analytic continuation in the lower
${\Pi}^{-}$ and the upper ${\Pi}^{+}$ half-plane, respectively, and $\Lambda(x)$ is a diagonal matrix
$$
\Lambda(x) = diag \left\{\left(\frac{x - i}{x + i}\right)^{{\ae}_1}, \ldots, \left(\frac{x - i}{x + i}\right)^{{\ae}_n}\right\},\; {\ae}_1, \ldots, {\ae}_n\in {\mathbb Z}.
$$

\noindent (6) Partial indices ${\ae}_1 \geq \ldots \geq {\ae}_n$ are stable, i.e.
\begin{equation}
\label{cond6}
{\ae}_1 - {\ae}_n \leq 1.
\end{equation}

The matrices of the following form constitute a simple subclass of class ${\mathcal S} K_2$:
\begin{equation}
\label{example}
G_{\varphi}(x) = \left(
\begin{array}{lr}
p(x) & q(x) e^{i \varphi x} \\
 q(x) e^{-i \varphi x} & p(x)
 \end{array}
 \right).
\end{equation}
In this {particular} case
\begin{equation}
\label{Ex1_6}
R_{\varphi}(x) =
\left(
\begin{array}{cc}
e^{\frac{i \varphi x}{2}} & 0 \\
& \\
0 & e^{-\frac{i \varphi x}{2}}
\end{array}
\right),\quad {F(x) =
\left(
\begin{array}{cc}
p(x) & q(x) \\
& \\
q(x) & p(x)
\end{array}
\right).}
\end{equation}
This function appears after Fourier transforms of the Wiener-Hopf equation, describing a problem in fracture mechanics.

We note that the matrix-functions of this type do not belong to any known class
of matrix-functions which admit explicit factorization (see, e.g. \cite{MishRog15a} and references therein).


\section{An algorithm}

\subsection{General construction}
\label{general}

By assumption any matrix $G_{\varphi}(x)\in {\mathcal S} K_n$ can be written in the form
\begin{equation}
\label{constr1}
G_{\varphi}(x) = F^{-}(x) \left(F^{-}(x)\right)^{-1} G_{\varphi} \left(F^{+}(x)\right)^{-1} F^{+}(x) =: F^{-}(x) G_{1, \varphi} F^{+}(x),
\end{equation}
where $F^{-}(x), F^{+}(x)$ are components of the factorization of the corresponding matrix $F$ (see (\ref{cond5})), and   the matrix $G_{1, \varphi}(x)$ is represented in the form (see (\ref{cond2}), (\ref{cond4}))
$$
G_{1, \varphi}(x) = \left(F^{-}(x)\right)^{-1} G_{\varphi} \left(F^{+}(x)\right)^{-1} = \left(F^{-}(x)\right)^{-1} R_{\varphi} F R_{\varphi}^{-1} \left(F^{+}(x)\right)^{-1}.
$$

{In addition to (1)--(6) we assume.}

\noindent (7) There exist a small parameter $\varepsilon=\varepsilon(\varphi)$ (more exactly its value will be described later), such that for all $x\in {\mathbb R}$ and any finite $\varphi$
\begin{equation}
\label{constr2}
G_{1, \varphi}(x) = \Lambda(x) + \varepsilon N_{\varphi}(x),
\end{equation}
with $\Lambda(x)$ determined by (\ref{cond5}) and matrix $N_{\varphi}(x)$
being bounded and locally H\"older continuous on ${\mathbb R}$.

Note that by assumption each entry of the matrix $N_{\varphi}(x)$ has a limit {as} $|x| \rightarrow +\infty$, i.e. there exists the value $N_{\varphi}(\infty) = 0$, i.e. $N_{\varphi}(x)\in H_{\mu}^{0}(\overline{\mathbb R})$.\footnote{As we can see later, the elements of our asymptotic factorization formula do not necessarily retain this property.}
Note also that the commutativity of the involved matrices is {not} assumed.

{Let us} look for the first step in asymptotic factorization of the matrix $G_{1, \varphi}(x)$  in the form
\begin{equation}
\label{repr1}
G_{1, \varphi}(x) = \left(I + \varepsilon N_{1,\varphi}^{-}(x) (\Lambda^{+}(x))^{-1}\right) \Lambda(x) \left(I + \varepsilon (\Lambda^{-}(x))^{-1} N_{1,\varphi}^{+}(x)\right),
\end{equation}
where
\begin{equation}
\label{repr2}
\Lambda(x) = \Lambda^{+}(x) \Lambda^{-}(x),
\end{equation}
with
\begin{equation}
\label{repr3}
\Lambda^{\pm}(x) = diag \left\{\left(\frac{x - i}{x + i}\right)^{{\ae}_1^{\pm}}, \ldots, \left(\frac{x - i}{x + i}\right)^{{\ae}_n^{\pm}}\right\},
\end{equation}
where
\begin{equation}
\label{repr4}
{\ae}_j^{+} = \max \{0, {\ae}_j\},\;\;\; {\ae}_j^{-} = \min \{0, {\ae}_j\},\;\; j=1,\ldots, n.
\end{equation}

In stable case, when condition (\ref{cond6}) holds, we have two possibilities, namely, ${\ae}_j \geq 0, \forall j = 1,\ldots, n$, or ${\ae}_j \leq 0, \forall j = 1,\ldots, n$.

$1^{0}$. Let us consider first the special case of matrix functions in ${\mathcal S} K_n$ for which the partial indices of factorization (\ref{cond5}) are such that ${\ae}_j = 1, \forall j = 1,\ldots, k$, ${\ae}_j = 0, \forall j = k+1,\ldots, n$, $1 \leq k \leq n$. Denote by $\Lambda_{0}$ the corresponding diagonal matrix in (\ref{cond5})
\begin{equation}
\label{alg1}
\Lambda_{0}(x) = diag \left\{\underbrace{\left(\frac{x - i}{x + i}\right), \ldots, \left(\frac{x - i}{x + i}\right)}\limits_{k\; {\mathrm{times}}}, 1, \ldots, 1\right\}.
\end{equation}

In this case
\begin{equation}
\label{alg2}
\Lambda^{+}_{0}(x) = \Lambda_{0}(x),\;\;\; \Lambda^{-}_{0}(x) = I.
\end{equation}
In this notation representation (\ref{repr1}) is equivalent to
\begin{equation}
\label{repr5}
G_{1, \varphi}(x) = \left(\Lambda^{+}_{0}(x) + \varepsilon N_{1,\varphi}^{-}(x)\right) \left(I + \varepsilon  N_{1,\varphi}^{+}(x)\right).
\end{equation}

Comparing terms for the lower powers of $\varepsilon$ we get the following matrix boundary value problem for determination of factors
$N_{1,\varphi}^{-}$, $N_{1,\varphi}^{+}$:
\begin{equation}
\label{bvp1}
\Lambda^{+}_{0}(x) N_{1,\varphi}^{+}(x) + N_{1,\varphi}^{-}(x) =  N_{\varphi}(x), \; x\in {\mathbb R}.
\end{equation}
It is customary to denote
$$
{M}_{0,\varphi}(x) \equiv N_{\varphi}(x).
$$
Since the matrix coefficient $\Lambda^{+}_{0}$ of problem (\ref{bvp1}) is a diagonal matrix, then this problem is equivalent to $n^2$ independent scalar boundary value problems
\begin{equation}
\label{bvp2}
\left\{
\begin{array}{lll}
\left(\frac{x - i}{x + i}\right) n_{11}^{+}(x) + n_{11}^{-}(x)&=&m_{0,11}, \\
&\ldots& \\
\left(\frac{x - i}{x + i}\right) n_{1n}^{+}(x) + n_{1n}^{-}(x)&=&m_{0,1n}, \\
&\ldots& \\
\left(\frac{x - i}{x + i}\right) n_{k1}^{+}(x) + n_{k1}^{-}(x)&=&m_{0,k1}, \\
&\ldots& \\
\left(\frac{x - i}{x + i}\right) n_{kn}^{+}(x) + n_{kn}^{-}(x)&=&m_{0,kn}, \\
n_{k+1;1}^{+}(x) + n_{k+1;1}^{-}(x)&=&m_{0,k+1;1}, \\
&\ldots& \\
n_{k+1;n}^{+}(x) + n_{k+1;n}^{-}(x)&=&m_{0,k+1;n}, \\
&\ldots& \\
n_{n1}^{+}(x) + n_{n1}^{-}(x)&=&m_{0,n1}, \\
&\ldots& \\
n_{nn}^{+}(x) + n_{nn}^{-}(x)&=&m_{0,nn},
\end{array}
\right.
\end{equation}
where $n_{pq}^{+}, n_{pq}^{-}, m_{0,pq}$ are entries of matrices $N_{1,\varphi}^{+}$, $N_{1,\varphi}^{-}$, ${M}_{0,\varphi}$, respectively.

Let us introduce a collection of pairs of analytic functions
\begin{equation}
\label{an_pq}
\Omega^{\pm}_{0}[m_{0,pq}](z) := \frac{z - i}{2 \pi i} \int\limits_{-\infty}^{\infty} \frac{m_{0,pq}(t) d t}{(t - i)(t - z)},\;\;\; z\in \Pi^{\pm}\;\;\; (1 \leq p, q \leq n).
\end{equation}
One can directly check that the introduced functions $\Omega^{\pm}_{0}[m_{0,pq}](z)$ are analytic in the respective half-plane, and satisfy the identity
\begin{equation}
\label{standard_fac}
\Omega^{+}_{0}[m_{0,pq}](x)-\Omega^{-}_{0}[m_{0,pq}](x)=m_{0,pq}(x), \quad x\in {\mathbb R},
\end{equation}
and the following additional condition holds true:
\begin{equation}
\label{cond_i}
\Omega^{+}_{0}[m_{0,pq}](i)=0.
\end{equation}
Using integrals (\ref{an_pq}) we get representation of bounded solutions to boundary value problems (\ref{bvp2}).

{
Bounded
solutions to the first $k\times n$ boundary value problems can be delivered by the formulas (see, e.g. \cite[p. 120]{Gak77})
\begin{equation}
\label{sol1_1}
n_{pq}^{+}(z) = \frac{z + i}{z - i} \Omega^{+}_{0}[m_{0,pq}](z),\quad n_{pq}^{-}(z) = - \Omega^{-}_{0}[m_{0,pq}](z),\;
\end{equation}
for $1 \leq p \leq k, 1 \leq q \leq n$,
and their boundary values are represented as
\begin{equation}
\label{sol1_3}
n_{pq}^{+}(x) = \frac{1}{2} \frac{x + i}{x - i} \Big(m_{0,pq}(x) + S_{0}[m_{0,pq}](x)\Big),
\end{equation}
\[
n_{pq}^{-}(x) = \frac{1}{2} \Big(m_{0,pq}(x) - S_{0}[m_{0,pq}](x)\Big),
\]
where $S_{0}$ is a corrected singular integral operator on the real line (see, e.g. \cite[p. 51-52]{Gak77})\footnote{We use the same notation for the Cauchy type operator $\Omega^{\pm}_{0}$ and the singular integral operator $S_{0}$ applied to matrices: $\Omega^{\pm}_{0}[M_{0,\varphi}](z)$, $S_{0}[M_{0,\varphi}](x)$.}
\begin{equation}
\label{sing_op}
S_{0}[m_{0,pq}](x) = \frac{x - i}{\pi i} \int\limits_{-\infty}^{\infty} \frac{m_{0,pq}(t) dt}{(t - i)(t - x)}.
\end{equation}
}

{On the other hand, bounded solutions to the remaining boundary value problems in  (\ref{bvp2}) can be delivered by the formulas
\begin{equation}
\label{sol1_2}
n_{pq}^{\pm}(z) = \pm \Omega^{\pm}_{0}[m_{0,pq}](z),\quad k + 1 \leq p \leq n, 1 \leq q \leq n.
\end{equation}
and their boundary values are consequently represented as
\begin{equation}
\label{sol1_4}
n_{pq}^{\pm}(x) = \frac{1}{2} \Big(m_{0,pq}(x) \pm S_{0}[m_{0,pq}](x)\Big).
\end{equation}}


It is known (see, e.g. \cite[p. 51-52]{Gak77}) that operator $S_{0}$ is a bounded operator on the space $H_{\mu}(\overline{\mathbb R})$ of bounded H\"older continuous functions on the open real line ${\mathbb R}$. The norm of this operator will be denoted\footnote{For the standard singular integral operator on the real line, i.e. the Hilbert transform, an exact values of its norm in H\"older spaces are known, see \cite{Ale75}.}
\begin{equation}
\label{sing_op_norm}
C_{\mu} := \|S_{0}\|_{H_{\mu}(\overline{\mathbb R}) \rightarrow H_{\mu}(\overline{\mathbb R})}.
\end{equation}

Let us refine the factorization of the matrix $G_{1, \varphi}(x)$, i.e. look for {a presentation of $G_{1, \varphi}(x)$} in the form
\begin{equation}
\label{repr2_1}
G_{1, \varphi}(x) = \left(I + \varepsilon N_{1,\varphi}^{-}(x) (\Lambda^{+}_{0}(x))^{-1} + \varepsilon^2 N_{2,\varphi}^{-}(x) (\Lambda^{+}_{0}(x))^{-1}\right) \times
\end{equation}
$$
\times \Lambda_{0}(x) \left(I + \varepsilon N_{1,\varphi}^{+}(x) + \varepsilon^2 N_{2,\varphi}^{+}(x)\right),
$$
where $N_{1,\varphi}^{-}(x)$, $N_{1,\varphi}^{+}(x)$ are those found at the previous step.

As before, relation (\ref{repr2_1}) is equivalent to
\begin{equation}
\label{repr2_2}
\Lambda_{0}(x) + \varepsilon N_{\varphi}(x) = \left(\Lambda^{+}_{0}(x) + \varepsilon N_{1,\varphi}^{-}(x) + \varepsilon^2 N_{2,\varphi}^{-}(x)\right)
\left(I + \varepsilon N_{1,\varphi}^{+}(x) + \varepsilon^2 N_{2,\varphi}^{+}(x)\right).
\end{equation}

Comparing terms at $\varepsilon^0, \varepsilon^1, \varepsilon^2$ we get  the following matrix boundary value problem for determination of factors
$N_{2,\varphi}^{-}$, $N_{2,\varphi}^{+}$:
\begin{equation}
\label{bvp_2}
 \Lambda^{+}_{0}(x) N_{2,\varphi}^{+}(x) + N_{2,\varphi}^{-}(x) = {{M}_{1,\varphi}(x)\equiv - N_{1,\varphi}^{-}(x) N_{1,\varphi}^{+}(x),}\; x\in {\mathbb R},
\end{equation}
{where the  right hand-side is already known}. {The solution to} this problem is given by the formulas similar to (\ref{sol1_1}) and (\ref{sol1_2}),
in which $m_{0,pq}$ should be replaced by $m_{1,pq}$. Essential difference {in comparison to the previous step} is that the functions $m_{1,pq}(x)$ is no longer vanishing at infinity, but remain to be bounded. Anyway, the functions $\Omega^{\pm}_{0}[m_{1,pq}](z)$ are well defined, i.e. analytic in the respective domains, having H\"older continuous boundary values, and satisfy (\ref{standard_fac}),
(\ref{cond_i}).


One can proceed in the same manner. Thus on the $r$-th step we use the representation
$$
G_{1, \varphi}(x) = \Lambda_{0}(x) + \varepsilon N_{\varphi}(x) =
$$
\begin{equation}
\label{repr_r}
= \left(I + \varepsilon N_{1,\varphi}^{-}(x) (\Lambda^{+}_{0}(x))^{-1} + \ldots + \varepsilon^r N_{r,\varphi}^{-}(x) (\Lambda^{+}_{0}(x))^{-1}\right) \times
\end{equation}
$$
\times \Lambda_{0}(x) \left(I + \varepsilon N_{1,\varphi}^{+}(x) + \ldots + \varepsilon^r N_{r,\varphi}^{+}(x)\right),
$$
where $N_{1,\varphi}^{-}(x), \ldots, N_{r-1,\varphi}^{-}(x)$, $N_{1,\varphi}^{+}(x), \ldots, N_{r-1,\varphi}^{+}(x)$ are found {in} the previous steps.
{This} leads to the matrix boundary value problem
\begin{equation}
\label{bvp_r}
\Lambda^{+}_{0}(x)  N_{r,\varphi}^{+}(x) + N_{r,\varphi}^{-}(x) = {M}_{r-1,\varphi}, \; x\in {\mathbb R},
\end{equation}
where
$$
{M}_{r-1,\varphi} = -\left[N_{1,\varphi}^{-}(x) N_{r-1,\varphi}^{+}(x) +
N_{2,\varphi}^{-}(x) N_{r-2,\varphi}^{+}(x) + \ldots + N_{r-1,\varphi}^{-}(x) N_{1,\varphi}^{+}(x)\right].
$$
{The solution to} this problem is given by {a} formulas similar to (\ref{sol1_1}) and (\ref{sol1_2}).

{Finally,} the factorization of the matrix function $G_{1, \varphi}(x)$ is given in the form
\begin{equation}
\label{asy_fact_gen}
G_{1, \varphi}(x) = G^-_{1,\varepsilon}(x) \Lambda_0(x) G^+_{1,\varepsilon}(x),
\end{equation}
where the respective minus and plus matrix functions are presented {by their asymptotic series}
\begin{equation}
\label{asy_fact}
 G^-_{1,\varepsilon}(x)= {I}+\sum\limits_{r=1}^{\infty} \varepsilon^r N_{r,\varphi}^{-}(x)  (\Lambda^{+}_{0}(x))^{-1},\quad
 G^+_{1,\varepsilon}(x)={I}+\sum\limits_{r=1}^{\infty} \varepsilon^r N_{r,\varphi}^{+}(x).
\end{equation}
Here $N_{r,\varphi}^{-}(x)$, $N_{r,\varphi}^{+}(x)$ are bounded solutions to the matrix problems
(\ref{bvp_r}) for any $r\in {\mathbb N}$. {From} the properties of the solutions to (\ref{bvp_r}) it follows that
$N_{r,\varphi}^{-}(x)$, $N_{r,\varphi}^{+}(x)$ are rather {finite than vanishing at infinity  ($|x| \rightarrow +\infty$).}

\vspace{3mm}
$2^{0}$. In general case of stable indices we have ${\ae}_j = s + 1, \forall j = 1,\ldots, k$, ${\ae}_j = s, \forall j = k+1,\ldots, n$, $1 \leq k \leq n$, where $s\in {\mathbb Z}$ is any integer number. In this case
$$
\Lambda(x) = diag \left\{\underbrace{\left(\frac{x - i}{x + i}\right)^{s+1}, \ldots, \left(\frac{x - i}{x + i}\right)^{s+1}}\limits_{k\; {\mathrm{times}}}, \left(\frac{x - i}{x + i}\right)^{s}, \ldots, \left(\frac{x - i}{x + i}\right)^{s}\right\},
$$
hence
$$
\Lambda(x) = \left(\frac{x - i}{x + i}\right)^{s} \Lambda_0(x),
$$
where {$\Lambda_0(x)=\Lambda_0^+(x)$} is the diagonal matrix given by (\ref{alg1}).

Let us consider first an auxiliary relation
\begin{equation}
\label{gen1}
\Lambda_0(x) + \varepsilon \left(\frac{x + i}{x - i}\right)^s N_{\varphi}(x) = \left(I + \varepsilon \tilde{N}_{1,\varphi}^{-}(x) (\Lambda^{+}_{0}(x))^{-1}\right) \Lambda_{0}(x) \left(I + \varepsilon \tilde{N}_{1,\varphi}^{+}(x)\right).
\end{equation}
Then we arrive at the following matrix boundary value problem:
\begin{equation}
\label{bvp_gen1}
\Lambda_0(x) \tilde{N}_{1,\varphi}^{+}(x) + \tilde{N}_{1,\varphi}^{-}(x) {=\tilde{M}_{0,\varphi}(x)\equiv  \left(\frac{x + i}{x - i}\right)^s N_{\varphi}(x).}
\end{equation}

{Then, using the same line of the reasoning as for the previous case, one receives the
asymptotic formula for factorization of the matrix-function (\ref{constr2}) in the form
\begin{equation}
\label{asy_fact_gen2}
G_{1, \varphi}(x) = \tilde{G}^-_{1,\varepsilon}(x) \Lambda(x) \tilde{G}^+_{1,\varepsilon}(x),
\end{equation}
where the respective minus and plus factors are presented in forms of asymptotic series
\begin{equation}
\label{asy_fact2}
 \tilde{G}^-_{1,\varepsilon}(x)= {I}+\sum\limits_{r=1}^{\infty} \varepsilon^r \tilde{N}_{r,\varphi}^{-}(x)  (\Lambda^{+}_{0}(x))^{-1},\quad
 G^+_{1,\varepsilon}(x)={I}+\sum\limits_{r=1}^{\infty} \varepsilon^r \tilde{N}_{r,\varphi}^{+}(x).
\end{equation}}

\subsection{Convergence}
\label{convergence}

The following Theorem gives conditions when{, in the case of the stable partial indices, the proposed} asymptotic factorization becomes  an explicit one, i.e. gives convergence conditions for the asymptotic series involved.

\begin{theorem}
\label{converge}
{Let $G_{\varphi}$ be a matrix which meets the conditions (1) to (7). Let the parameter $\varepsilon$ (defined in (7))} satisfies the inequality
\begin{equation}
\label{conv_cond}
|\varepsilon| \leq 1/A
\end{equation}
with the constant {$A=A(\varphi)$} being equal to
\begin{equation}
\label{conv_cond1}
A = \|N_{\varphi}(\cdot)\|_{\mu} (1 + C_{\mu})^2,
\end{equation}
$\|N_{\varphi}(\cdot)\|_{\mu}$ being the norm of the matrix function $N_{\varphi}(x)$ in the H\"older space $H_{\mu}(\overline{\mathbb R})$,\footnote{This matrix norm can be any {\it sub-multiplicative} norm of
square matrix-functions in $H_{\mu}(\overline{\mathbb R})$, i.e. satisfying an inequality $\|B(\cdot) \cdot C(\cdot)\|_{\mu} \leq \|B(\cdot)\|_{\mu} \|C(\cdot)\|_{\mu}$ for any pair of square matrix-functions  $B(x), C(x)$
(e.g. any {\it induced operator norm}, see e.g. \cite{HorJoh85}).}
and $C_{\mu}$ being the norm of the singular integral operator ${S}_{0} : H_{\mu}(\overline{\mathbb R}) \rightarrow H_{\mu}(\overline{\mathbb R})$.

Then both series in  (\ref{asy_fact}) converge for all $x\in {\mathbb R}$.
\end{theorem}
$\triangleleft$
It suffices to prove Theorem only in the above considered case $1^{0}$ (see Subsec. \ref{general}).
It follows (see, e.g., \cite[p. 48]{Gak77}) that singular integral operator $S_{0}$ is bounded in H\"older spaces since the ``standard'' singular integral operator ${S}$ is. Let $C_{\mu}$ be the norm of $S_{0}$ in H\"older space $H_{\mu}(\overline{\mathbb R})$.
Then we have the following series of estimates
$$
\|N_{1,\varphi}^{\mp}(\cdot)\|_{\mu} \leq \alpha_1 \|N_{\varphi}(\cdot)\|_{\mu} (1 + C_{\mu}), \; {\mathrm{where}}\; \alpha_1 = \frac{1}{2},
$$
$$
\|N_{2,\varphi}^{\mp}(\cdot)\|_{\mu} \leq \frac{1}{2} \|M_{1,\varphi}(\cdot)\|_{\mu} (1 + C_{\mu}),
$$
and
$$
\|M_{1,\varphi}(\cdot)\|_{\mu} \leq \left(\alpha_1 \|N_{\varphi}(\cdot)\|_{\mu} (1 + C_{\mu})\right)^2,
$$
i.e.
$$
\|N_{2,\varphi}^{\mp}(\cdot)\|_{\mu} \leq \alpha_2 \|N_{\varphi}(\cdot)\|_{\mu}^{2} (1 + C_{\mu})^{3}, \; {\mathrm{where}}\; \alpha_2 = \frac{1}{2} \alpha_1^2.
$$
Finally, for each $r \geq 2$
$$
\|N_{r,\varphi}^{\mp}(\cdot)\|_{\mu} \leq \frac{1}{2} \|M_{r-1,\varphi}(\cdot)\|_{\mu} (1 + C_{\mu}),
$$
and
$$
\|M_{r-1,\varphi}(\cdot)\|_{\mu} \leq \left(\alpha_1 \alpha_{r-1} + \alpha_2 \alpha_{r-2} + \ldots + \alpha_{r-1} \alpha_1\right) \|N_{\varphi}(\cdot)\|_{\mu}^{r} \left(1 + C_{\mu}\right)^{2r - 2},
$$
i.e.
$$
\|N_{r,\varphi}^{\mp}(\cdot)\|_{\mu} \leq \alpha_r \|N_{\varphi}(\cdot)\|_{\mu}^{r} (1 + C_{\mu})^{2r - 1}, \; {\mathrm{where}}\; \alpha_r = \frac{1}{2} (\alpha_1 \alpha_{r-1} + \ldots + \alpha_{r-1} \alpha_1).
$$
Few first coefficients $\alpha_r$ we can calculate explicitly, namely, $\alpha_1 = 1/2$, $\alpha_2 = 1/8$, $\alpha_3 = 1/16$. As for coefficients with large enough
indices we can proof by induction that
$$
\alpha_r < \frac{1}{16(r-3)}, \;\;\; \forall r \geq 12.
$$
Therefore
$$
\|\varepsilon^r N_{r,\varphi}^{\mp}(\cdot)\|_{\mu} \leq |\varepsilon|^r \frac{1}{16(r-3)(1 + C_{\mu})} \left(\|N_{\varphi}(\cdot)\|_{\mu} (1 + C_{\mu})^2\right)^{r}, \; \forall r \geq 12.
$$
Since the sequence $\sqrt[r]{\frac{1}{16(r-3)(1 + C_{\mu})}} \leq 1$ is increasing for sufficiently large $r$ and
$$
\lim\limits_{r \rightarrow\infty} \sqrt[r]{\frac{1}{16(r-3)(1 + C_{\mu})}} = 1,
$$
then (since $\|(\Lambda^{+}(x))^{-1}\| = 1$) the convergence of the series (\ref{asy_fact})
for all $x\in {\mathbb R}$ follows from (\ref{conv_cond}).
$\triangleright$

\begin{remark}
\label{rem_conv}
It follows from the standard properties of the Cauchy type integral and singular integral with Cauchy kernel
that conditions of Theorem \ref{converge} guarantee convergence of the series in the right-hand side of (\ref{asy_fact})
in the half-planes ${\Pi}^{-}$, ${\Pi}^{+}$, respectively.
\end{remark}

\begin{remark}
\label{infinity}
In fact, the decay of the second term in the right-hand side of (\ref{constr2}) at infinity
follows from the properties of matrices
of the considered class and the proposed construction. Meanwhile, the solutions of matrix boundary
problems (\ref{bvp_r}) do not retain such behavior at infinity.\footnote{In contrast to the case of canonical factorization in
\cite{MishRog14a}, where the solutions of the corresponding problems at each step of factorization vanish at infinity.}
\end{remark}

\begin{remark}
\label{smallness}
If the number $A = A(\varphi)$ in Theorem \ref{converge} is small enough, i.e.
\begin{equation}
\label{conv_cond1_without}
A = \|N_{\varphi}(\cdot)\|_{\mu} (1 + C_{\mu})^2<1,
\end{equation}
then the results remains valid for $\varepsilon = 1$ and the described procedure {will work} without any changes.
\end{remark}

\subsubsection{The behavior of boundary values of solution at infinity}

Let us recall (see \cite[$\S$ 4]{Gak77}) few facts on the behavior of the Cauchy type operator $\Omega_{0}$ and singular operator $S_{0}$, presenting these properties in the matrix form.


{1) Let $n\times n$ matrix function $M\in H_{\mu}^{0}(\overline{\mathbb R})$, i.e. is H\"older continuous on $\overline{\mathbb R}$ and vanishing at infinity. The the Cauchy type integral $\Omega_{0}^{\pm}[M](z)$ (see (\ref{an_pq}))
determines two analytic matrices having boundary functions $\Omega_{0}^{\pm}[M](x)\in H_{\mu}^{0}(\overline{\mathbb R})$ and satisfies the Sokhotsky-Plemelj formula
\begin{equation}
\label{inf2}
\Omega_{0}^{+}[M](x) - \Omega_{0}^{-}[M](x) = M(x), \;\;\; x\in {\mathbb R},
\end{equation}
and vanishes at $z=i$:
\begin{equation}
\label{inf3}
\Omega_{0}^{+}[M](i) = 0.
\end{equation}
}

Since
$$
\Omega_{0}^{\pm}[M](z) = \frac{1}{2\pi i} \int\limits_{-\infty}^{\infty} \frac{M(t) d t}{(t - z)} -
\frac{1}{2\pi i} \int\limits_{-\infty}^{\infty} \frac{M(t) d t}{(t - i)}
$$
then it follows from (\ref{inf2}) {that there exists}
\begin{equation}
\label{inf4}
 \lim\limits_{z \rightarrow \infty} \Omega_{0}^{\pm}[M](z) = - \Omega_{0}^{+}[M](i) = 0.
\end{equation}

The corresponding singular integral operator $S_{0}[M](x)$ (see (\ref{sing_op}) satisfy the relation
\begin{equation}
\label{inf7}
 \lim\limits_{x \rightarrow \infty} S_{0}[M](x) = - 2 \Omega_{0}^{+}[M](i) = 0.
\end{equation}
It follows from (\ref{inf4}) and another form of  the Sokhotsky-Plemelj formula
$$
2\Omega_{0}^{\pm}[M](x) = \pm  M(x) +  S_{0}[M](x).
$$

2) We can consider (\ref{inf2}) as the boundary value problem (jump problem). Then
any other bounded analytic solution $\tilde{\Omega}_{0}^{\pm}[M](z)$ to this problem (different from $\Omega_{0}^{\pm}[M](z)$) has the form
\begin{equation}
\label{inf5}
\tilde{\Omega}_{0}^{\pm}[M](z) = \Omega_{0}^{\pm}[M](z) +\tilde{\Omega}_{0}^{\pm}[M](i),
\end{equation}
i.e. differs of our special solution $\Omega_{0}^{\pm}[M](z)$ in a constant matrix.
The above considerations show {that there exist}
\begin{equation}
\label{inf4a}
 \lim\limits_{z \rightarrow \infty} \tilde{\Omega}_{0}^{\pm}[M](z) = - \tilde{\Omega}_{0}^{+}[M](i),
\end{equation}
\begin{equation}
\label{inf7a}
 \lim\limits_{z \rightarrow \infty} \tilde{S}_{0}[M](z) = - 2 \tilde{\Omega}_{0}^{+}[M](i).
\end{equation}

\begin{remark}
\label{inf_sol}
Any solution $N_{1,\varphi}^{+}, N_{1,\varphi}^{-}$ to the problem (\ref{bvp1}) has the following behavior at infinity
\begin{equation}
\label{inf_sol1}
\lim\limits_{z \rightarrow \infty} N_{1,\varphi}^{\pm}(z) = \mp \tilde{\Omega}_{0}^{+}[M_{0,\varphi}](i),
\end{equation}
where $\tilde{\Omega}_{0}^{\pm}[M_{0,\varphi}]$ is a chosen bounded solution to the jump problem
$$
\Omega^{+}(x) + \Omega^{-}(x) = M_{0,\varphi}(x), \;\;\; x\in {\mathbb R}.
$$
If, in particular, $\tilde{\Omega}_{0}^{\pm}[M_{0,\varphi}] = {\Omega}_{0}^{\pm}[M_{0,\varphi}]$,
then $N_{1,\varphi}^{+}, N_{1,\varphi}^{-}$ vanishes at infinity.

The matrix $M_{1,\varphi}(x) = - N_{1,\varphi}^{-}(x) N_{1,\varphi}^{+}(x)$ has the following asymptotics at infinity
\begin{equation}
\label{inf_sol2}
\lim\limits_{x \rightarrow \infty} M_{1,\varphi}(x) = \left[\tilde{\Omega}_{0}^{+}[M_{0,\varphi}](i)\right]^2,
\end{equation}
and is vanishing at infinity either matrix $\tilde{\Omega}_{0}^{+}[M_{0,\varphi}](i)$ or its square is the zero-matrix.
\end{remark}

\begin{remark}
\label{unique}
The algorithm described in Subsec. \ref{general} is not unique (cf. \cite{MishRog14a}) since solution to the problem (\ref{bvp_r})
at each step of approximation
is defined up to the constant matrix
\begin{equation}
\label{unique_r}
N_{r,\varphi}^{+} = \left[\Lambda^{+}\right]^{-1} \left\{\tilde{\Omega}_{r - 1}^{+}[M_{r - 1,\varphi}] - C_{r - 1}\right\},\;
N_{r,\varphi}^{-} = - \left\{\tilde{\Omega}_{r - 1}^{-}[M_{r - 1,\varphi}] - C_{r - 1}\right\},
\end{equation}
where $\tilde{\Omega}_{r - 1}^{\pm}[M_{r - 1,\varphi}]$ is one of possible solutions to the jump problem
$$
\Omega^{+}(x) + \Omega^{-}(x) = M_{r - 1,\varphi}(x), \;\;\; x\in {\mathbb R},
$$
and $C_{r - 1}$ is a constant matrix.

The first $k$ rows of are prescribed by equalities
$$
C_{r - 1,pq} = \tilde{\Omega}_{r - 1}^{+}[m_{r - 1,pq}](i), \;\;\; 1 \leq p \leq k,\;\; 1\leq q \leq n,
$$
but the remaining $n - k$ rows can be chosen arbitrary.

The choice made in Secs. \ref{general}, \ref{convergence} corresponds to
\begin{equation}
\label{unique_m}
C_{r - 1} = \tilde{\Omega}_{r - 1}^{+}[M_{r - 1,\varphi}](i), \;\;\; \forall r = 0, 1, 2 \ldots
\end{equation}
\end{remark}

{\begin{remark}
\label{asymp_proc}
According to Remark~\ref{unique},  is it possible to minimize the norm of
the matrices $M_{j,\varphi}(x)$, at each consequent step ($j=1,2,\ldots N$) adopting the arbitrary constants. As a results, one obtains a factorization
with required boundary behavior and asymptotic properties. The norms of the one or both factors $N^\pm_{j,\varphi}(x)$ can be controlled. However, this would not always
guarantee the improvement for the convergence as the
norm of matrix $M_{j,\varphi}(x)$ to be normalized on the next step may well increase as the result. This suggests various scenarios for the optimization. One can also choose the arbitrary constant to preserve some local property of the solution (for example, to set a particular value of some components of the factors $G^\pm(\varphi,x)$ at a points $x=x_s$ (depending on the number of arbitrary constants). Moreover, limiting values of the matrices at infinity ($x=\infty$) also can be controlled in this way.
All these fact allow to construct specific factorization with properties suitable for a particular application.
\end{remark}}

\begin{remark}
\label{some_relations}
Interestingly, whatever factorization is constructed,
the factors satisfy the following conditions
\begin{equation}
\label{as_factors}
g^-_{pq}(-i)=\delta_{pq},\quad 1\le p \le n,\;\; 1\le q\le k;\quad G^-_{1,\varepsilon}(\infty)\cdot G^+_{1,\varepsilon}(\infty)=I.
\end{equation}
Here $G^-_{1,\varepsilon}(x)=\left\{g^-_{pq}(x)\right\}$, $p,q=1,2,\ldots, n$ and the integer $k$ is defined in (\ref{alg1}).

{The authors believe that impossibility to set the uniqueness condition in the form $G^-_{1,\varepsilon}(\infty)=G^+_{1,\varepsilon}(\infty)=I$ (as it was in the algorithm for the canonical factorization \cite{MishRog14a}) is related to the fact that the proposed algorithms
always produces the conditions (\ref{as_factors})$_1$ automatically and thus it would impose too much restrictions to the factor if additionally $G^-_{1,\varepsilon}(\infty)=I$ is imposed. If any other algorithm for factorization of that class of the matrix functions is proposed which is free from the condition (\ref{as_factors})$_1$
then probably some kind of a uniqueness result would be available.}
And finally, if one restricts the procedure to $N$ asymptotic term only, the latter conditions will be valid with an accuracy of the order $O(\varphi^N)$.

All the above statements of this Remark are valid in the general case too (see (\ref{asy_fact2})) with the replacement of the matrices $G^\pm_{1,\varepsilon}(x)$ by $\tilde G^\pm_{1,\varepsilon}(x)$.
\end{remark}

\section{Example of a $2\times2$ nonrational matrix}
\label{Example}

\subsection{Asymptotic algorithm}

{\bf Example 1.} Let us consider the matrix function on the real line ${\mathbb R}$
\begin{equation}
\label{Ex1_1}
F(x) :=
\left(
\begin{array}{cc}
\frac{x - 3 i}{x + i} & \frac{- 2 i}{x + i} \\
& \\
\frac{4 i}{x + i} & \frac{x + 3 i}{x + i}
\end{array}
\right).
\end{equation}
$F(x)$ satisfies the properties (1)--(7) (Sec. \ref{class}), i.e. $F(x)\in {\mathcal S} K_2$ in particular,
 it admits the right factorization
$$
F(x) = F^{-}(x) \Lambda(x) F^{+}(x), \;\;\; x\in {\mathbb R},
$$
with
\begin{equation}
\label{Ex1_2}
F^{-}(x) =
\left(
\begin{array}{cc}
1 & - 1 \\
& \\
- 1 & 2
\end{array}
\right),\; F^{+}(x) =
\left(
\begin{array}{cc}
2 & 1 \\
& \\
 1 & 1
\end{array}
\right),\; \Lambda(x) = {\Lambda_{0}(x) =}
\left(
\begin{array}{cc}
\frac{x - i}{x + i} & 0 \\
& \\
0 & 1
\end{array}
\right).
\end{equation}}
Partial indices of this factorization are equal
\begin{equation}
\label{Ex1_3}
{\ae}_1 = 1,\;\;\; {\ae}_2 = 0,
\end{equation}
and thus they are stable ${\ae}_1 - {\ae}_2 = 1$.

Then with $R_{\varphi}(x)$ as in (\ref{Ex1_6})
\begin{equation}
\label{Ex1_7}
G_{\varphi}(x) := R_{\varphi}(x) F(x) R_{\varphi}^{-1}(x) =
\left(
\begin{array}{cc}
\frac{x - 3 i}{x + i} & \frac{- 2 i e^{i \varphi x}}{x + i} \\
& \\
\frac{4 i e^{- i \varphi x}}{x + i} & \frac{x + 3 i}{x + i}
\end{array}
\right).
\end{equation}

By construction {(see (\ref{constr1}))}
$$
G_{1, \varphi}(x) := \left[F^{-}(x)\right]^{-1} G_{\varphi}(x) \left[F^{+}(x)\right]^{-1}.
$$
Here $\left[F^{-}(x)\right]^{-1} = F^{+}(x), \left[F^{+}(x)\right]^{-1} = F^{-}(x)$.
Thus
\begin{equation}
\label{Ex1_8}
G_{1, \varphi}(x) =
\left(
\begin{array}{cc}
\frac{x - 9 i + 4 i e^{- i \varphi x} + 4 i e^{i \varphi x}}{x + i} & \frac{12 i - 4 i e^{- i \varphi x} - 8 i e^{i \varphi x}}{x + i} \\
& \\
\frac{- 6 i + 4 i e^{- i \varphi x} + 2 i e^{i \varphi x}}{x + i} & \frac{x + 9 i - 4 i e^{- i \varphi x} - 4 i e^{i \varphi x}}{x + i}
\end{array}
\right).
\end{equation}
As in general case, the matrix function $G_{1, \varphi}(x)$ can be represented in the form
$$
G_{1, \varphi}(x) = \Lambda_{0}(x) + \varepsilon N_{\varphi}.
$$
We can represent the matrix $\varepsilon N_{\varphi}$ in the form involving trigonometric functions:
\begin{equation}
\label{Ex1_11}
\varepsilon N_{\varphi}(x) =
\left(
\begin{array}{cc}
\frac{- 16 i \sin^2 \frac{\varphi x}{2}}{x + i} & \frac{8 i \sin \frac{\varphi x}{2} \left(3 \sin \frac{\varphi x}{2} - i \cos \frac{\varphi x}{2}\right)}{x + i} \\
& \\
\frac{- 4 i  \sin \frac{\varphi x}{2} \left(3 \sin \frac{\varphi x}{2} + i \cos \frac{\varphi x}{2}\right)}{x + i} & \frac{16 i \sin^2\frac{\varphi x}{2}}{x + i}
\end{array}
\right).
\end{equation}
It follows from (\ref{Ex1_11}) that the matrix function $\varepsilon N_{\varphi}$ is ``small'', i.e. vanishing at infinity with properly chosen $\varepsilon$.

Following Remark \ref{smallness}, we do not separate the existing small parameter $\varepsilon$ from the matrices
but simply work the the entire matrices as they appear, i.e. consider representation 
$$
G_{1, \varphi}(x) = \Lambda_{0}(x) + M_{0,\varphi}(x).
$$
Now we give another representation of  $M_{0,\varphi}$ which is calculated directly
\begin{equation}
\label{Ex1_8a}
M_{0, \varphi}(x) =
\left(
\begin{array}{cc}
\frac{-8 i + 4 i e^{- i \varphi x} + 4 i e^{i \varphi x}}{x + i} & \frac{12 i - 4 i e^{- i \varphi x} - 8 i e^{i \varphi x}}{x + i} \\
& \\
\frac{- 6 i + 4 i e^{- i \varphi x} + 2 i e^{i \varphi x}}{x + i} & \frac{8 i - 4 i e^{- i \varphi x} - 4 i e^{i \varphi x}}{x + i}
\end{array}
\right).
\end{equation}
It can be factorized in the form
$$
M_{0, \varphi}(x) = M_{0, \varphi}^{+}(x) + M_{0, \varphi}^{-}(x),
$$
where
\begin{equation}
\label{Ex1_8+}
M_{0, \varphi}^{+}(x) =
2 i \left(
\begin{array}{cc}
\frac{-4 + 2 e^{i \varphi x} + 2 e^{- \varphi}}{x + i} & \frac{6 - 4 e^{i \varphi x} - 2 e^{- \varphi}}{x + i} \\
& \\
\frac{- 3 + e^{i \varphi x} + 2 e^{- \varphi}}{x + i} & \frac{4 - 2 e^{i \varphi x} - 2 e^{- \varphi}}{x + i}
\end{array}
\right),
\end{equation}
\begin{equation}
\label{Ex1_8-}
M_{0, \varphi}^{-}(x) =
2 i \left(
\begin{array}{cc}
\frac{2 (e^{- i \varphi x} - e^{- \varphi})}{x + i} & \frac{- 2 (e^{- i \varphi x} - e^{- \varphi})}{x + i} \\
& \\
\frac{2 (e^{- i \varphi x} - e^{- \varphi})}{x + i} & \frac{- 2 (e^{- i \varphi x} - e^{- \varphi})}{x + i}
\end{array}
\right).
\end{equation}

Let us find solution to the problem
\begin{equation}
\label{Ex1_13}
\Lambda^{+}_{0}(x) N_{1,\varphi}^{+}(x) + N_{1,\varphi}^{-}(x) =  M_{0, \varphi}(x),
\end{equation}
or what is equivalent to 4 separate scalar problems
\begin{equation}
\label{Ex1_14}
\left\{
\begin{array}{lll}
\frac{x - i}{x + i} n_{11}^{+} + n_{11}^{-}&=&m_{0,11}, \\
\frac{x - i}{x + i} n_{12}^{+} + n_{12}^{-}&=&m_{0,12}, \\
n_{21}^{+} + n_{21}^{-}&=&m_{0,21}, \\
n_{22}^{+} + n_{22}^{-}&=&m_{0,22},
\end{array}
\right.
\end{equation}
where $n_{kl}^{+}, n_{kl}^{-}, m_{0,kl}, k, l = 1,2,$ are entries of the matrices  $ N_{1,\varphi}^{+}, N_{1,\varphi}^{-}, M_{0, \varphi}$, respectively. We look for (partial) solution to the boundary value problems (\ref{Ex1_14}) in the class of bounded analytic functions.

The bounded solutions to the problems (\ref{Ex1_14}) can be written in the form
$$
\begin{array}{ll}
n_{11}^{+}(z) = \frac{z + i}{z - i} \left(m_{0,11}^{+}(z) - c_{11}\right), & n_{11}^{-}(z) = m_{0,11}^{-} + c_{11}, \\
n_{12}^{+}(z) = \frac{z + i}{z - i} \left(m_{0,12}^{+}(z) - c_{12}\right), & n_{12}^{-}(z) = m_{0,12}^{-} + c_{12}, \\
n_{21}^{+}(z) = m_{0,21}^{+}(z) - c_{21}, & n_{21}^{-}(z) = m_{0,21}^{-} + c_{21}, \\
n_{22}^{+}(z) = m_{0,22}^{+}(z) - c_{22}, & n_{22}^{-}(z) = m_{0,22}^{-} + c_{22},
\end{array}
$$
where $c_{kl}, k,l = 1, 2,$ are appropriate constants which serve to make $n_{kl}^{+}(z)$ to be analytic in $\Pi^{+}$.

Note that since {the matrix $M_{0,\varphi}^{\pm}(x)$ tends to zero as $x \rightarrow \infty$, by construction, } then, for any choice {of the constant } matrix $C_{0}$, {the matrix from the right hand side of (\ref{bvp_2}), i.e.}
\begin{equation}
\label{ex_m1}
M_{1,\varphi}(x) = - N_{1,\varphi}^{-}(x) N_{1,\varphi}^{+}(x)
\end{equation}
has the following behavior at infinity
\begin{equation}
\label{ex_inf}
M_{1,\varphi}(\infty)\equiv \lim\limits_{x \rightarrow \infty} M_{1,\varphi}(x) = {C_{0}^2.}
\end{equation}}

We suggested before (see Remark \ref{unique}) to chose constant matrix $C_{0} = \{c_{0,pq}\}_{p,q = 1, 2}$ in such a way that
$$
C_{0} = M_{0, \varphi}^{+}(i).
$$
It gives {($C_{0}^{(0)}=C_{0}$)}
\begin{equation}
\label{Matrix_C_0}
{ C_{0}^{(0)}} =\Psi(\varphi)
\begin{pmatrix}
4 & -6 \cr
3 & -4
\end{pmatrix},\quad \Psi(\varphi)= e^{- \varphi} - 1.
\end{equation}
As the result, the limiting value of the new matrix-function $M_{1\varphi}(x)$ (see (\ref{inf_sol2})) to be factorized on the next step (\ref{ex_inf}) is computed:
Note that
\begin{equation}
\label{Matrix_M_11}
M^{(0)}_{1,\phi}(\infty)=- 2 \Psi^2(\varphi) I.
\end{equation}

Note that according to the Remark \ref{unique}, the constants $c_{21}, c_{22}$ can be taken arbitrary. Just for illustration, we consider {here the following three other choices of the matrix $C_0$: 
\begin{equation}
\label{Matrix_C_12}
C_{0}^{(1)} = \Psi(\varphi)
\begin{pmatrix}
4 & -6 \cr
0 & 0
\end{pmatrix},\quad
C_{0}^{(2)} = \Psi(\varphi)
\begin{pmatrix}
4 & -6 \cr
0 & 4
\end{pmatrix}, \quad  C_{0}^{(3)} = \Psi(\varphi)
\begin{pmatrix}
4 & -6 \cr
8/3 & -4
\end{pmatrix}.
\end{equation}
Then, the corresponding limiting values of the respective matrix functions $M_{1\varphi}(x)$ are
\begin{equation}
\label{Matrix_M_12}
M^{(1)}_{1,\phi}(\infty)= {8\Psi^2(\varphi)}
\begin{pmatrix}
2 &-3\cr
0 & 0
\end{pmatrix},\quad
{M^{(2)}_{1,\phi}(\infty)=16\Psi^2(\varphi)
\begin{pmatrix}
1 &-3\cr
0 & 1
\end{pmatrix},\quad M^{(3)}_{1,\phi}(\infty)=0.}
\end{equation}

{We chose to consider those cases for the following reasons: the first one looks the simplest natural choice (all arbitrary constants to set zero value). The second one was selected as it preserves the trigonal matrix form (with the equal elements on the main diagonal) for the matrices $N_{1, \varphi}^{\pm}(x)$ and thus for all consecutive matrices in the procedure. Finally, the third case guarantees that the
matrix function $M^{(3)}_{1,\phi}(x)$ preserves the same property at infinity as {the initial matrix function} $M_{0,\phi}(x)$.
}

The respective first order approximation of the factorization problem (\ref{Ex1_13}) has the form}
\begin{equation}
\label{Ex1_8+fin+012}
N_{1, \varphi; j}^{+}(x) =(\Lambda_0^+)^{-1} \big({M_{0, \varphi}^{+}(x)}-C_{0}^{(j)}\big), \quad N_{1, \varphi; j}^{-}(x) = {M_{0, \varphi}^{-}(x)}+C_{0}^{(j)},
\end{equation}

In fact, we define here four different variants of the asymptotic expansion
$N_{1, \varphi}^{\pm}(x)=N_{1, \varphi; j}^{\pm}(x)$.
Any of these solutions gives the first order approximates for the factorization in the form
\begin{equation}
\label{Ex1_12a}
\Lambda_{0}(x) + M_{0,\varphi}(x) = H_{1,\varphi}^-(x) \Lambda_{0}(x) H_{1,\varphi}^+(x) + O(\varphi^2),\quad \varphi\to0,
 \end{equation}
where
\begin{equation}
\label{Ex1_12aa}
H_{1,\varphi}^-(x)= I + N_{1,\varphi}^{-}(x) (\Lambda_{0}^{+}(x))^{-1},\quad
H_{1,\varphi}^+(x)= I + N_{1,\varphi}^{+}(x),
\end{equation}
and the estimate $O(\varphi^2)$ is uniform with respect to the variable $x$.

For the next step of approximation we use representation similar to (\ref{repr2_2}) and  arrive at the following matrix boundary value problem
\begin{equation}
\label{Ex1_21}
\Lambda^{+}_{0}(x) N_{2,\varphi}^{+}(x) + N_{2,\varphi}^{-}(x) =  M_{1, \varphi}(x) \equiv  - N_{1,\varphi}^{-}(x) N_{1,\varphi}^{+}(x).
\end{equation}
Entries of the matrix $M_{1, \varphi}$ are calculated directly.

{By using solution to problem (\ref{Ex1_21}) we get approximation at the second step.
\begin{equation}
\label{Ex1_20a}
\Lambda_{0}(x) + N_{\varphi}(x) =H_{2,\varphi}^-(x)\Lambda_{0}(x)H_{2,\varphi}^+(x)+O(\varphi^3),\quad \varphi\to 0,
\end{equation}
\begin{equation}
\label{Ex1_20aa}
H_{2,\varphi}^-(x)= I + N_{1,\varphi}^{-}(x) (\Lambda_{0}^{+}(x))^{-1}+N_{2,\varphi}^{-}(x)(\Lambda_{0}^{+}(x))^{-1},
\end{equation}
$$
H_{2,\varphi}^+(x)= I + N_{1,\varphi}^{+}(x)+N_{2,\varphi}^{+}(x).
$$}
{The procedure can be naturally continued to an arbitrary order of the small parameter $\varphi$.}

\subsection{Numerical simulations}

\subsubsection{Estimate of the reminder: standard choice of the constant $C_0=M_{0,\varphi}(i)$}

First, we consider the case with the choice (\ref{Matrix_C_0}) of the constants, $C_0^{(0)}$, which guarantees the solvability of the problem (\ref{bvp1}) and was instrumental  to prove the convergency result. {Having the first asymptotic terms of the factorization procedure in form (\ref{Ex1_12aa}), one can proceed to obtain
the next consequent approximations (\ref{Ex1_20aa}) by factorization of the matrix functions
$M_{1,\varphi}(x)$. In this specific example, one can write the analytical representation of $M_{1,\varphi}(x)$ from (\ref{Ex1_21}) which,
however, become cumbersome enough. For this reason we present numerical results for the components of the matrix $M_{1,\varphi}(x)$ which highlight
the fact that they are of the order $O(\varphi^2)$ and uniformly bounded with respect to the
variable $x\in R$. On the other hand, the results can be treated as a measure of the factorization error if one decides to stop the asymptotic procedure
immediately after the first step and to consider (\ref{Ex1_12aa}) as the sought-for approximate solution.}

\begin{figure}[h!]
  \begin{center}
    \includegraphics [scale=0.45]{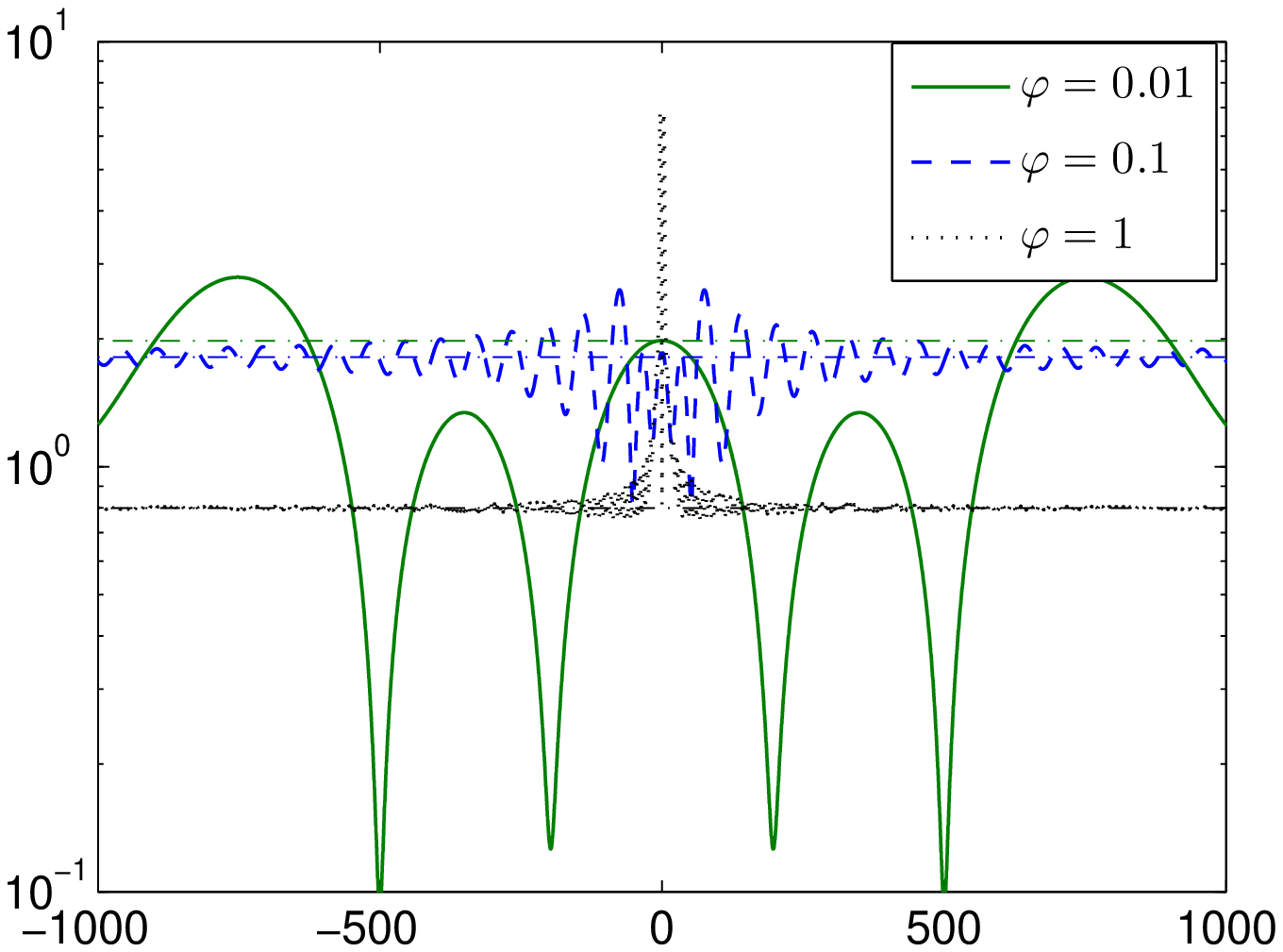}\includegraphics [scale=0.45]{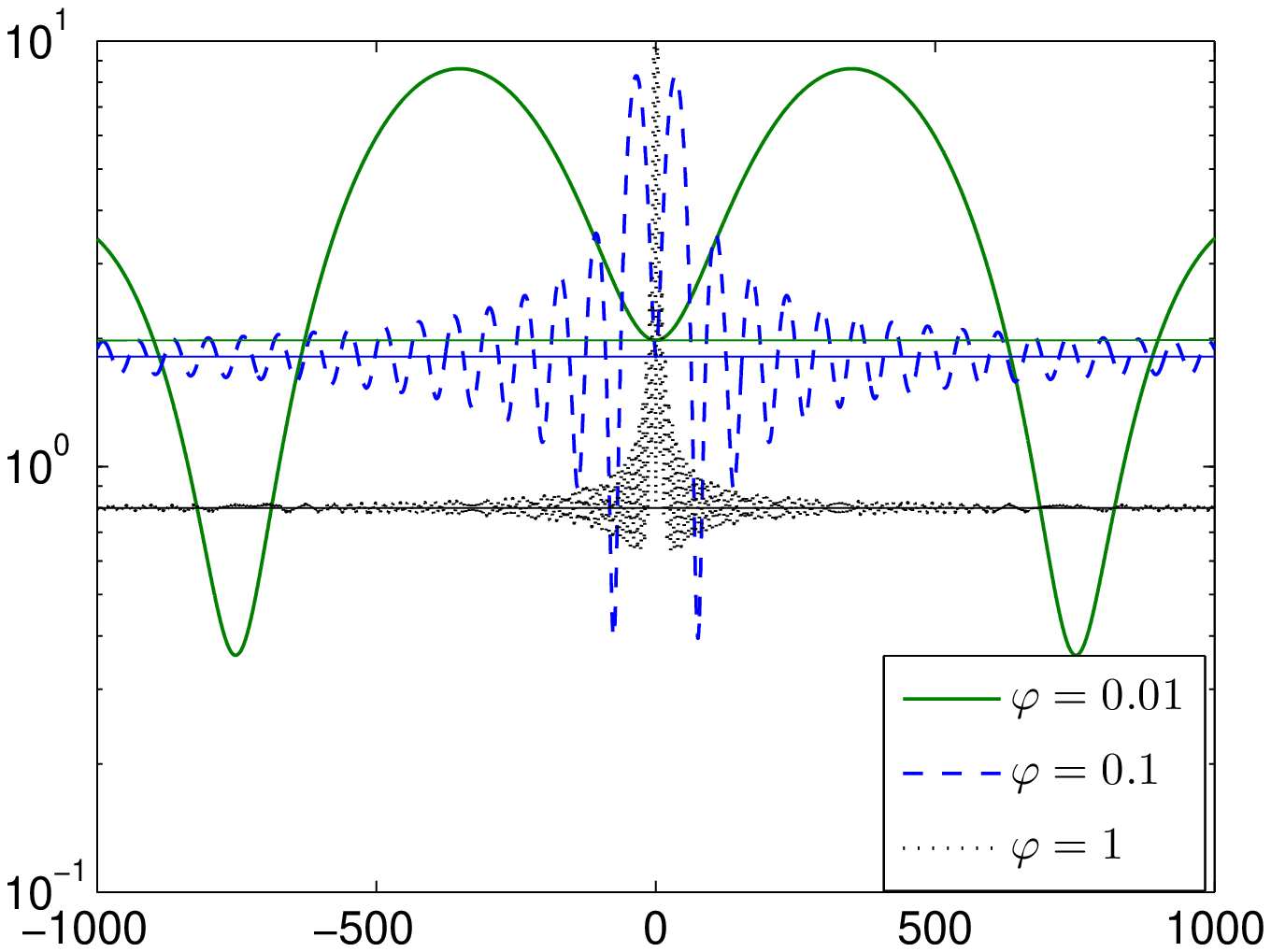}
       \put(-100,-5){\small$x$}
       \put(-365,115){$\varphi^{-2}m_{11}$}
       \put(-160,30){$\varphi^{-2}m_{22}$}
       \put(-300,-5){\small$x$}

  \end{center}
    \caption{\small Diagonal elements, $m_{jj}=m_{jj}(x,\varphi)$, $j=1,2$, of the matrix $M_{1,\varphi}(x)$ defined in (\ref{ex_m1}) for various values of the parameter $\varphi$.
    The elements are normalised to the value of the parameter $\varphi^2$. The horisontal lines show the limiting values of the normalised
    components at infinity $2(e^{-\varphi}-1)^2/\varphi^2$.}
\label{f1}
\end{figure}

In Fig.~\ref{f1} we present results concerning the two diagonal components, $m_{jj}(x,\varphi)$, $j=1,2$, of the matrix function $M_{1,\varphi}(x)$. Three different values of the
parameter $\varphi=1,0.1$ and $0.01$ are considered. Since, we expect that
\begin{equation}
\label{fi^2}
M_{1,\varphi}(x)=O(\varphi^2), \quad  \varphi\to0,
\end{equation}
the normalized values $\varphi^{-2}m_{jj}$ are given
instead. Note that the estimate $M_{1,\varphi}(\infty)=O(\varphi^2)$ follows immediately from (\ref{inf_sol2}) and the choice of the first row
in the matrix $C_0^{(0)}$.

\begin{figure}[h!]
  \begin{center}
    \includegraphics [scale=0.45]{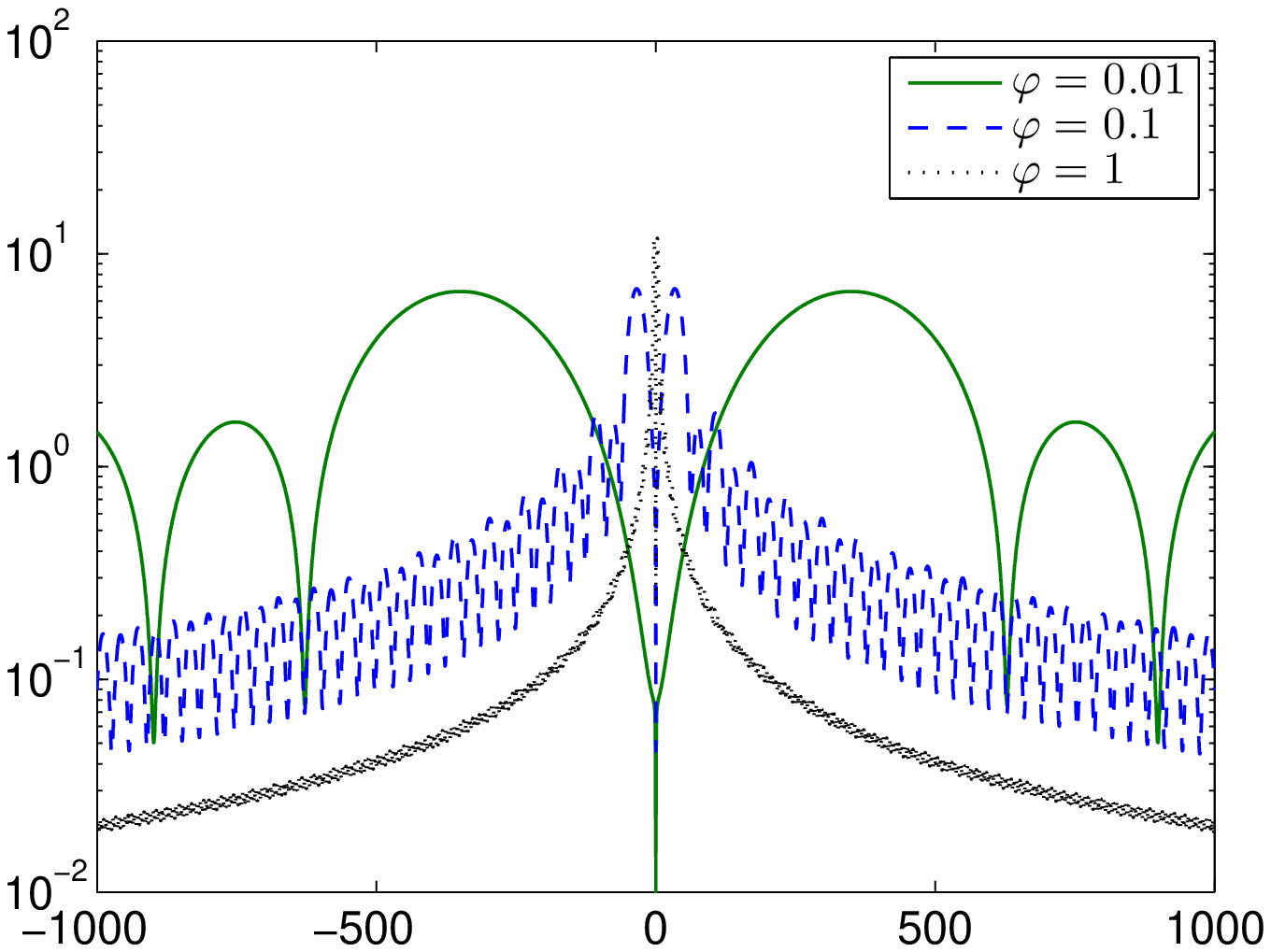}\includegraphics [scale=0.45]{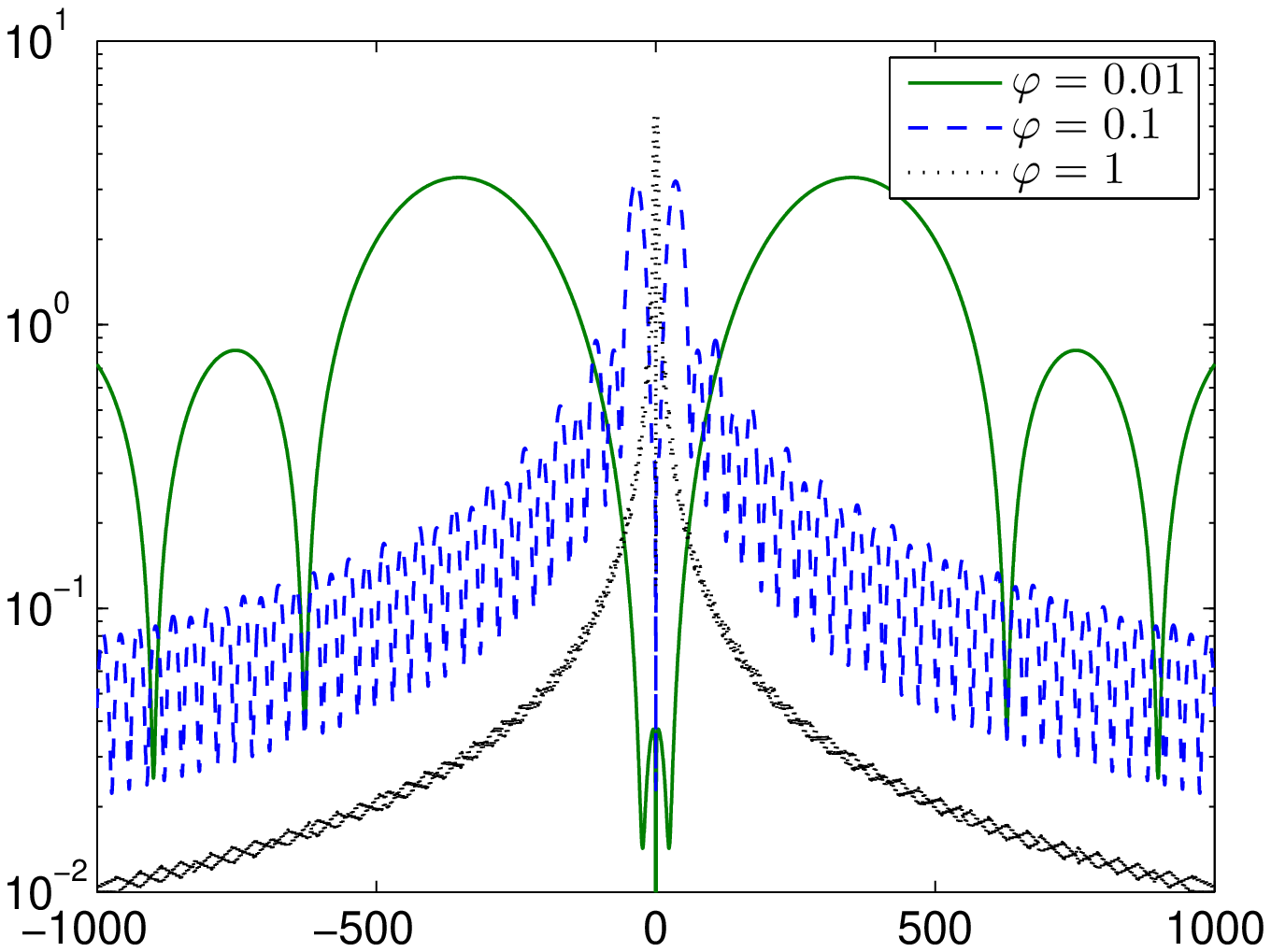}
       \put(-100,-5){\small$x$}
       \put(-370,115){$\varphi^{-2}m_{12}$}
       \put(-170,115){$\varphi^{-2}m_{21}$}
       \put(-300,-5){\small$x$}

  \end{center}
    \caption{\small The other two elements, $m_{ij}=m_{ij}(x,\varphi)$, $i+j=3$, of the matrix $M_{1,\varphi}(x)$ for different values of the parameter $\varphi=1;0.1;0.01$.
    It is clear that the both entries vanish at infinity $m_{ij}\to0$ as $|x|\to\infty$.}
\label{f2}
\end{figure}

Although the difference between the norms, $\|m_{jj}\|_\infty$, $j=1,2$, and the limiting values of the components, $|m_{jj}(\infty)|$, differ less with the decrease of the parameter $\varphi$ (see Fig.~\ref{f1}),  the limiting value of the matrix function, $M_{1,\varphi}(\infty)$, does not fully represent its max norm. For example,  $\|m_{22}\|_\infty> $$\|m_{11}\|_\infty$ apart from the fact that $|m_{22}(\infty)|=|m_{11}(\infty)|$.

\begin{figure}[h!]
  \begin{center}
    \includegraphics [scale=0.45]{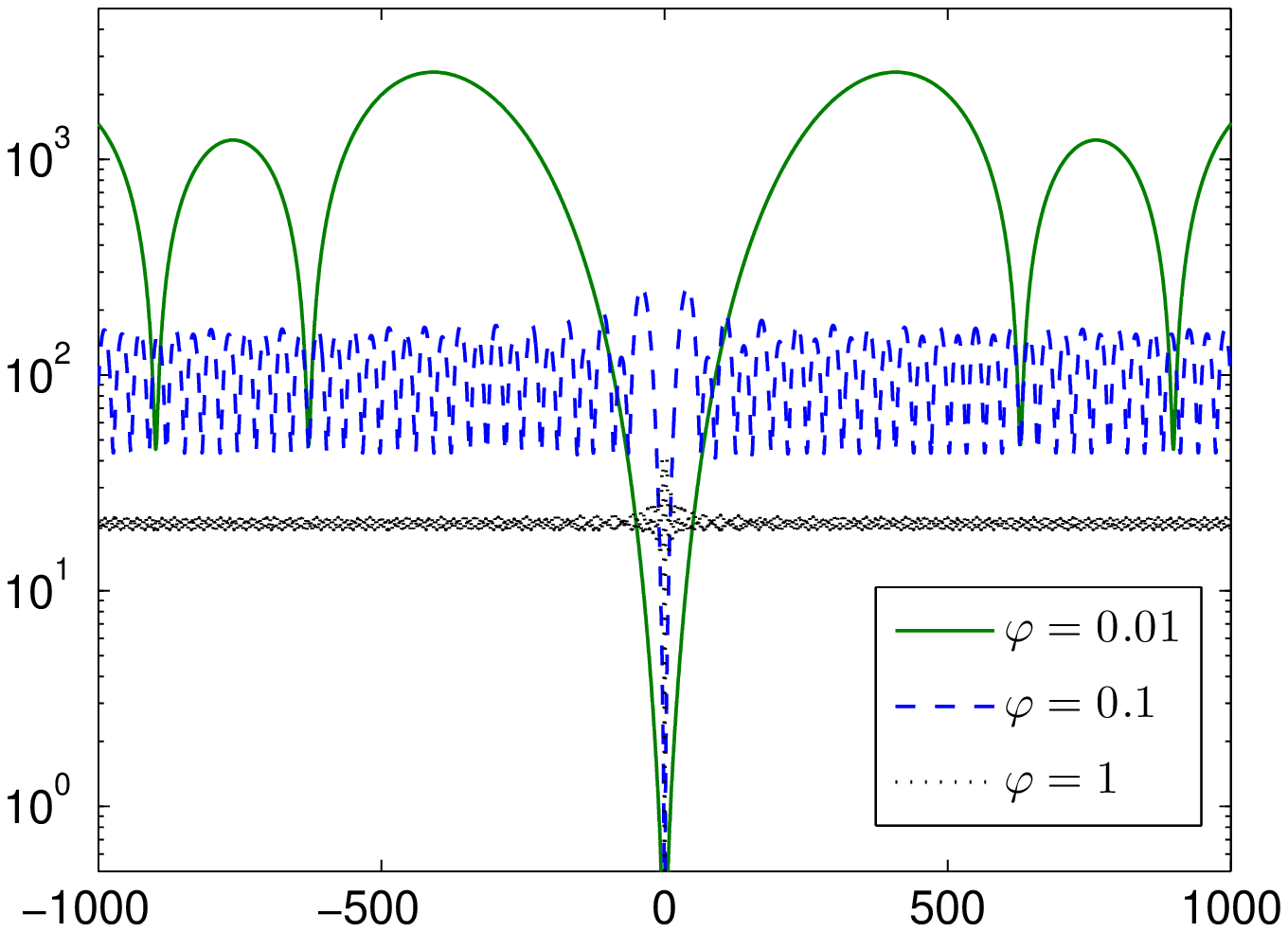}\includegraphics [scale=0.45]{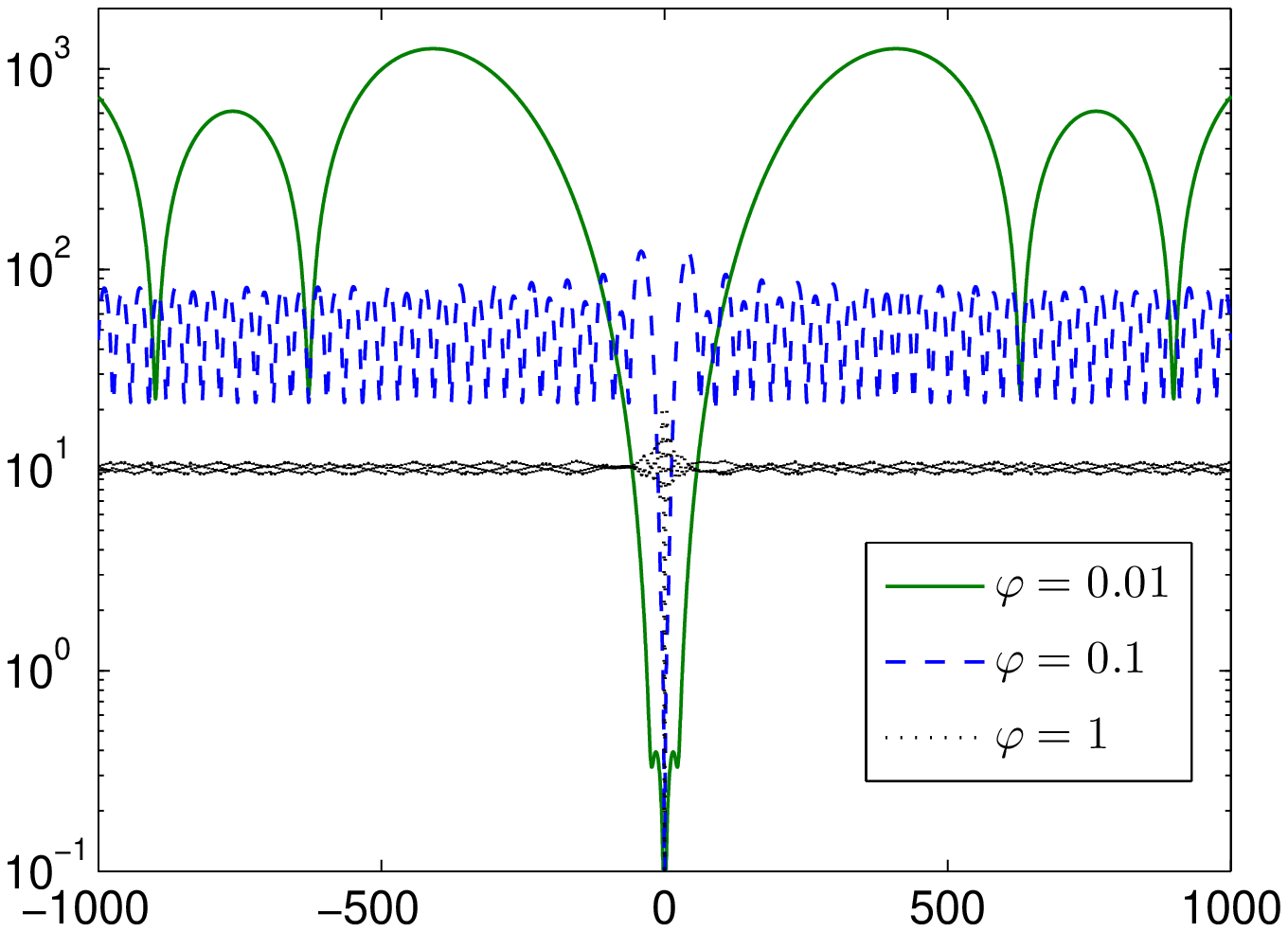}
       \put(-100,-5){\small$x$}
       \put(-370,40){$\displaystyle\frac{x}{\varphi^2}\,m_{12}(x)$}
       \put(-170,40){$\displaystyle\frac{x}{\varphi^2}\,m_{21}(x)$}
       \put(-300,-5){\small$x$}

  \end{center}
    \caption{\small Non diagonal elements, $m_{ij}(x)$, , $i+j=3$, of the matrix $M_{1,\varphi}(x)$ for various values of the parameter $\varphi$.
    The elements are normalised differently to demonstrate their asymptotics at infinity.
    }
\label{f3}
\end{figure}

The other two component, $m_{ij}$, $i+j=3$, of the matric function $M_{1,\varphi}(\infty)$ are presented in Fig.~\ref{f2}. To highlight the proven asymptotic estimate, here we aslo As expected, they
vanish at infinity as $x\to\infty$. One can also observe that $\|m_{12}\|_\infty> $$\|m_{21}\|_\infty$. Since, in the considered example, the functions are smooth at infinity (meromorfic function?) one can expect that
\begin{equation}
\label{contra}
m_{ij}(x,\varphi)=O(1/x),\quad x\to\infty.
\end{equation}
To demonstrate this fact, we also present in Fig.~\ref{f3} the normalised values of these components,  $x\varphi^{-2}\cdot m_{ij}(\varphi,x)$, $i+j=3$. Indeed the resulting functions are bounded and do not decay at infinity.

Surprisingly, their max norms grow while the value of the parameter $\varphi$ decreases and, from the first glance, this contradicts to (\ref{fi^2}) and (\ref{contra}).
In fact, it is all true as the components $m_{12}(\varphi,x)$ and $m_{21}(\varphi,x)$  behave like $x^{-1}\sin^2(\varphi x)$ and thus both estimates (\ref{fi^2}) and (\ref{contra}) are valid uniformly with respect to one or the other independent variable.

\subsubsection{Estimate of the reminder for the second choice of constant $C_0=C_0^{(1)}$}

\begin{figure}[h!]
  \begin{center}
    \includegraphics [scale=0.45]{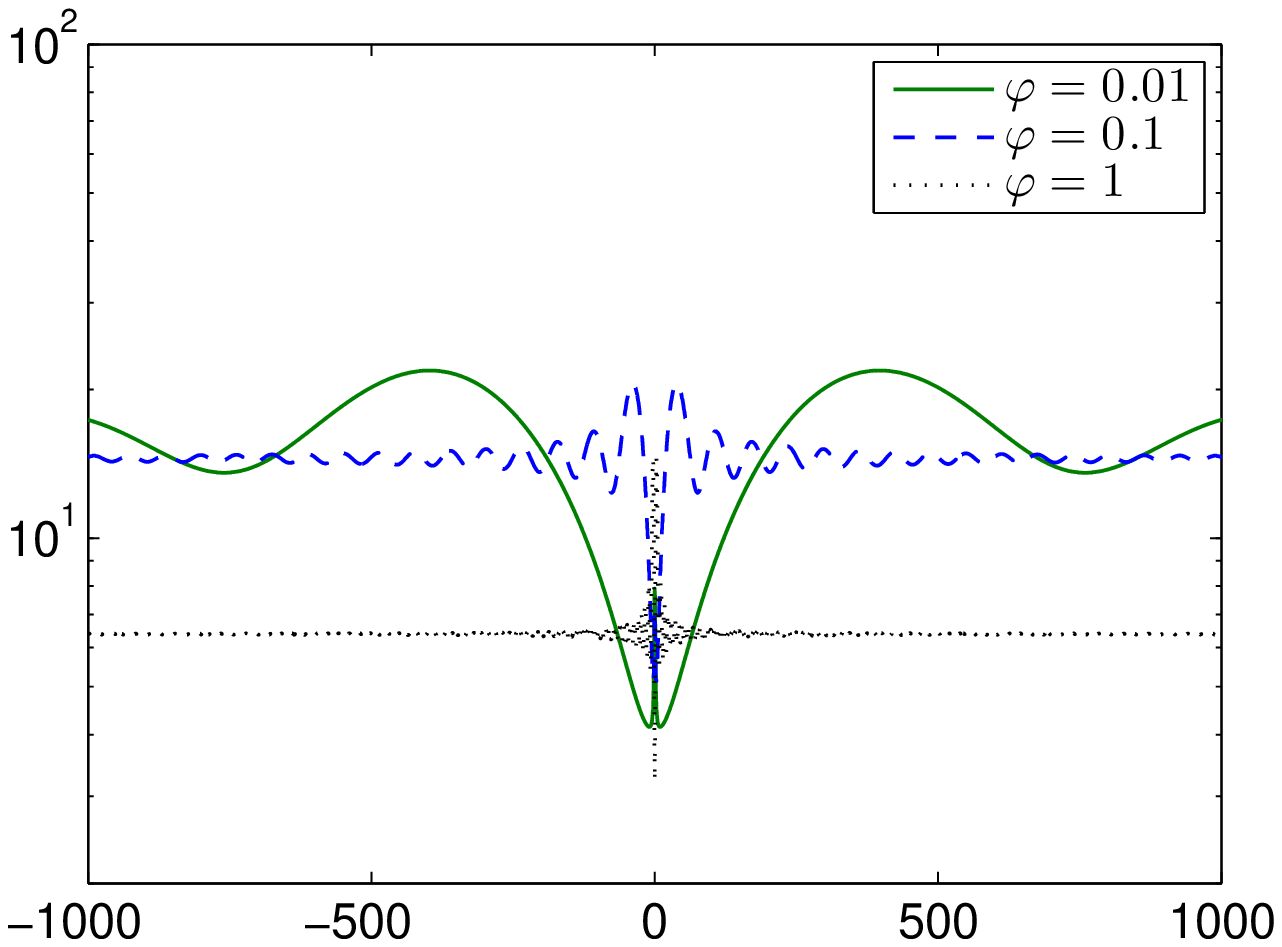}\includegraphics [scale=0.45]{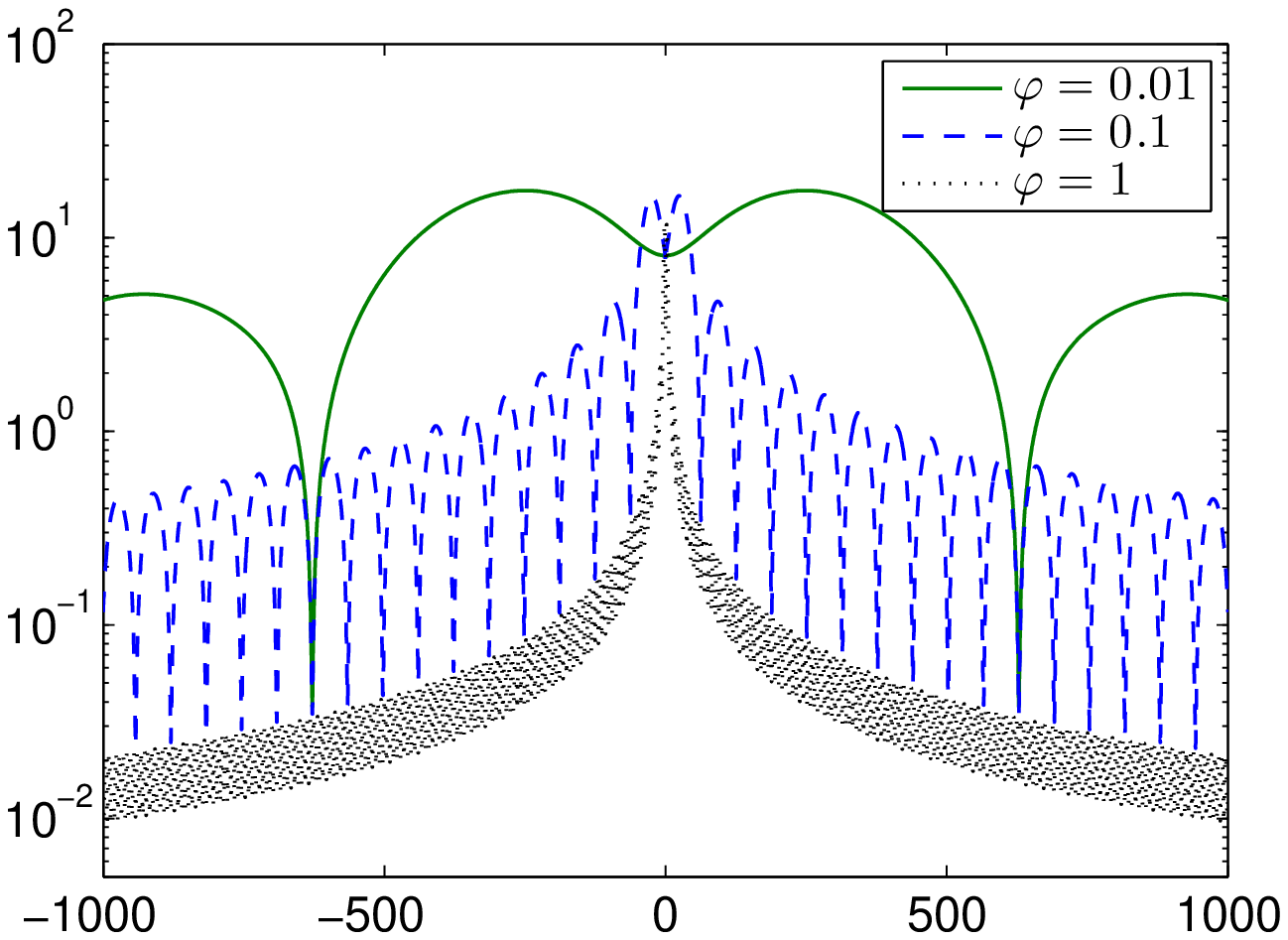}
       \put(-95,-5){\small$x$}
       \put(-350,110){$\varphi^{-2}m_{11}$}
       \put(-155,110){$\varphi^{-2}m_{22}$}
       \put(-285,-5){\small$x$}

  \end{center}
    \caption{\small Diagonal elements $m_{jj}(x)$ of the reminder matrix $M_{1,\varphi}(x)$ to be factorized on the next step for various values of the parameter $\varphi$. The elements are normalised to the value of the
    parameter $\varphi^2$. }
\label{f4}
\end{figure}

\begin{figure}[h!]
  \begin{center}
    \includegraphics [scale=0.45]{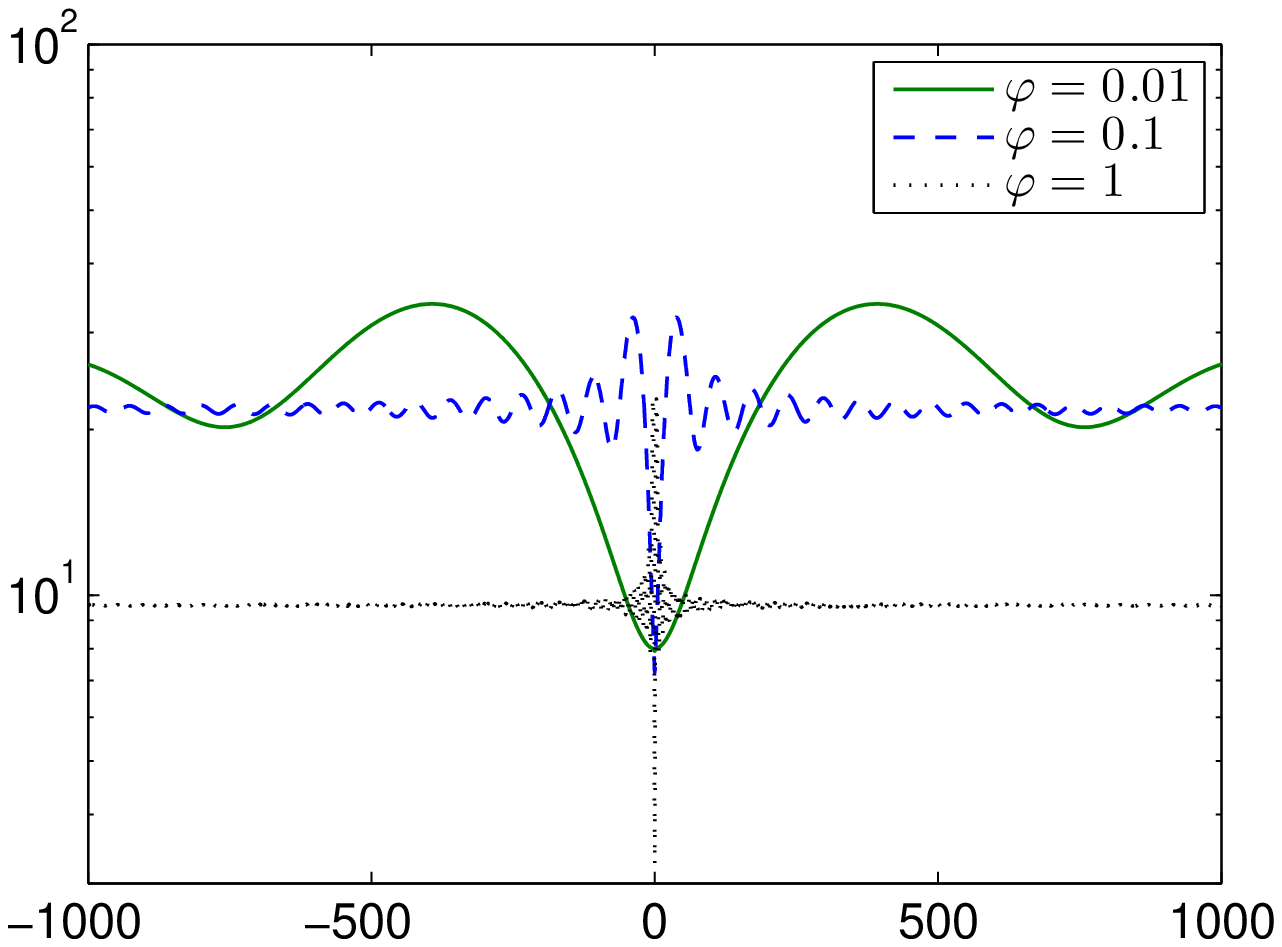}\includegraphics [scale=0.45]{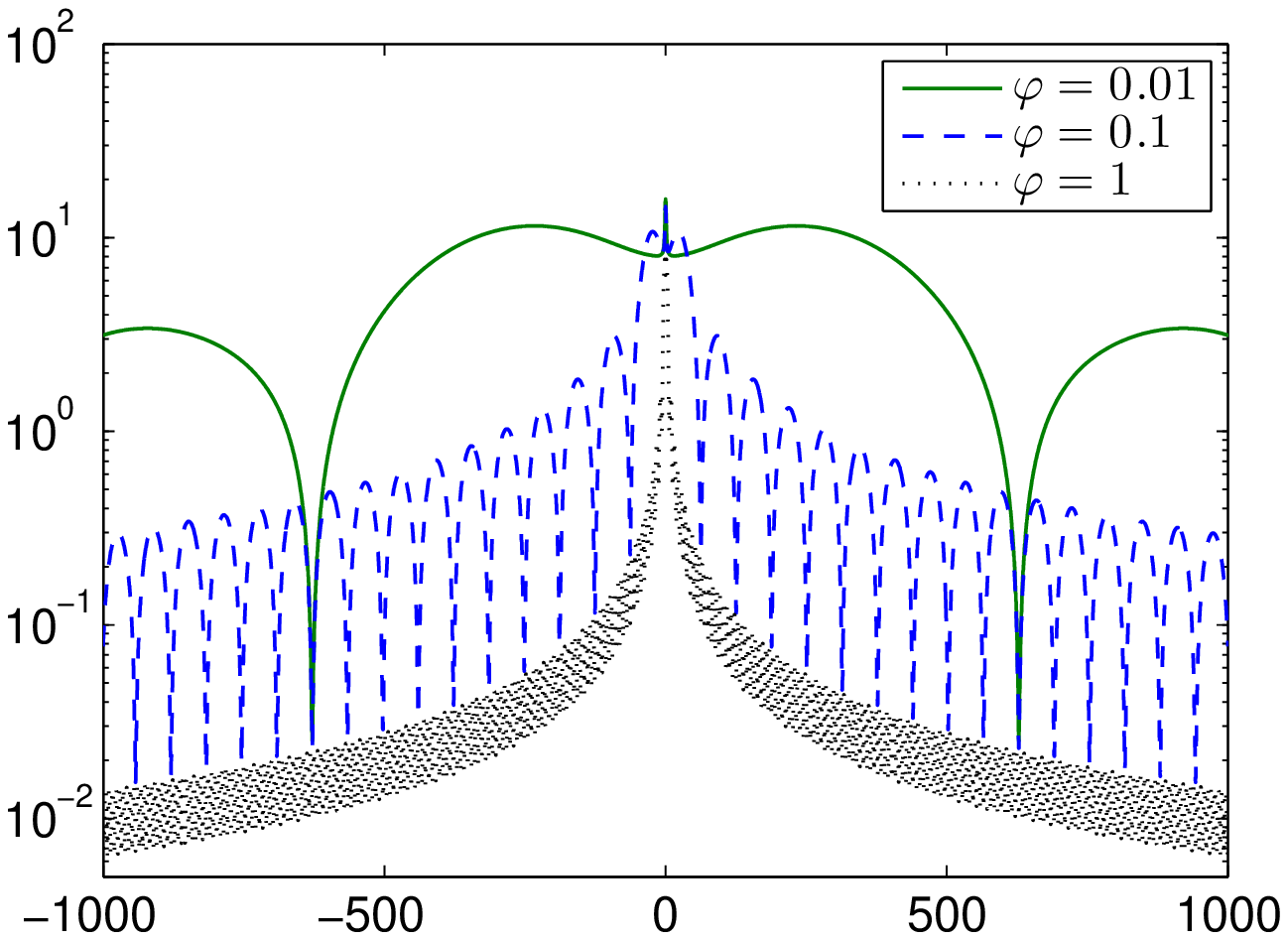}
       \put(-95,-5){\small$x$}
       \put(-350,110){$\varphi^{-2}m_{12}$}
       \put(-155,110){$\varphi^{-2}m_{21}$}
       \put(-285,-5){\small$x$}

  \end{center}
    \caption{\small The other two elements $m_{ij}(x)$ of the reminder matrix $M_{1,\varphi}(x)$ for various values of the parameter $\varphi$.
    The elements are normalised to the value of the parameter $\varphi^2$. }
\label{f5}
\end{figure}

Now we analyse the effect of a different choice of the constant matrix $C_0$ on the first step of the asymptotic procedure. Namely, we take it according to (\ref{Matrix_C_12})$_1$ which
provides the second simplest option when we set the arbitrary constants $c_{21}, c_{22}$ equal to zero.
In Fig.~\ref{f4} and Fig.~\ref{f5} present similar results as in the Fig.~\ref{f1} and Fig.~\ref{f2} above.

On the other hand, due to the structure of the chosen constant matrix $C_0^{(1)}$, the behaviour of the components in the first and the second lines are the same and their values are close to each other.

A a direct consequence of this choice, the components
$m_{21}(c)$, $m_{22}(x)$ of the matrix $M_{1,\varphi}(x)$ vanish at infinity ($x\to\infty$) and behave similarly to the component $m_{12}(x)$ and $m_{21}(x)$ from the previous example (compare (\ref{contra})). It is also obvious from Fig.~\ref{f4} and Fig.~\ref{f5} that the estimate (\ref{fi^2}) is valid.
However, and this was also suggested by the limiting value of the matrix $M_{1,\varphi}(\infty)$
in (\ref{Matrix_M_12}), the max norm of the matrix function $M_{1,\varphi}(x)$ increase approximately in ten times in comparison with the same value from the previous example.
Interestingly, the norm of the factors $N^\pm_{1,\varphi}(x)$ in both those cases are the same.

\subsubsection{Estimate of the reminder for the third choice of constant $C_0=C_0^{(2)}$}
It is important sometimes to have, if only possible, the factors well structured at infinity and, here, a triangular matrix with equal diagonal elements comes as a desirable choice.
This may be arranged by setting the constant matrix $C_0$ in the asymptotic algorithm by the relation (\ref{Matrix_C_12})$_2$, which guarantees solvability of the Hilbert-Riemann problem and, simultaneously, offers the both factors $N^\pm_{1,\varphi}(x)$ and the reminder $M_{1,\varphi}(x)$ represented by that type of the matrices at infinity (compare (\ref{Matrix_C_12})$_2$ and (\ref{Matrix_M_12}$_2$)).

\begin{figure}[h!]
  \begin{center}
    \includegraphics [scale=0.50]{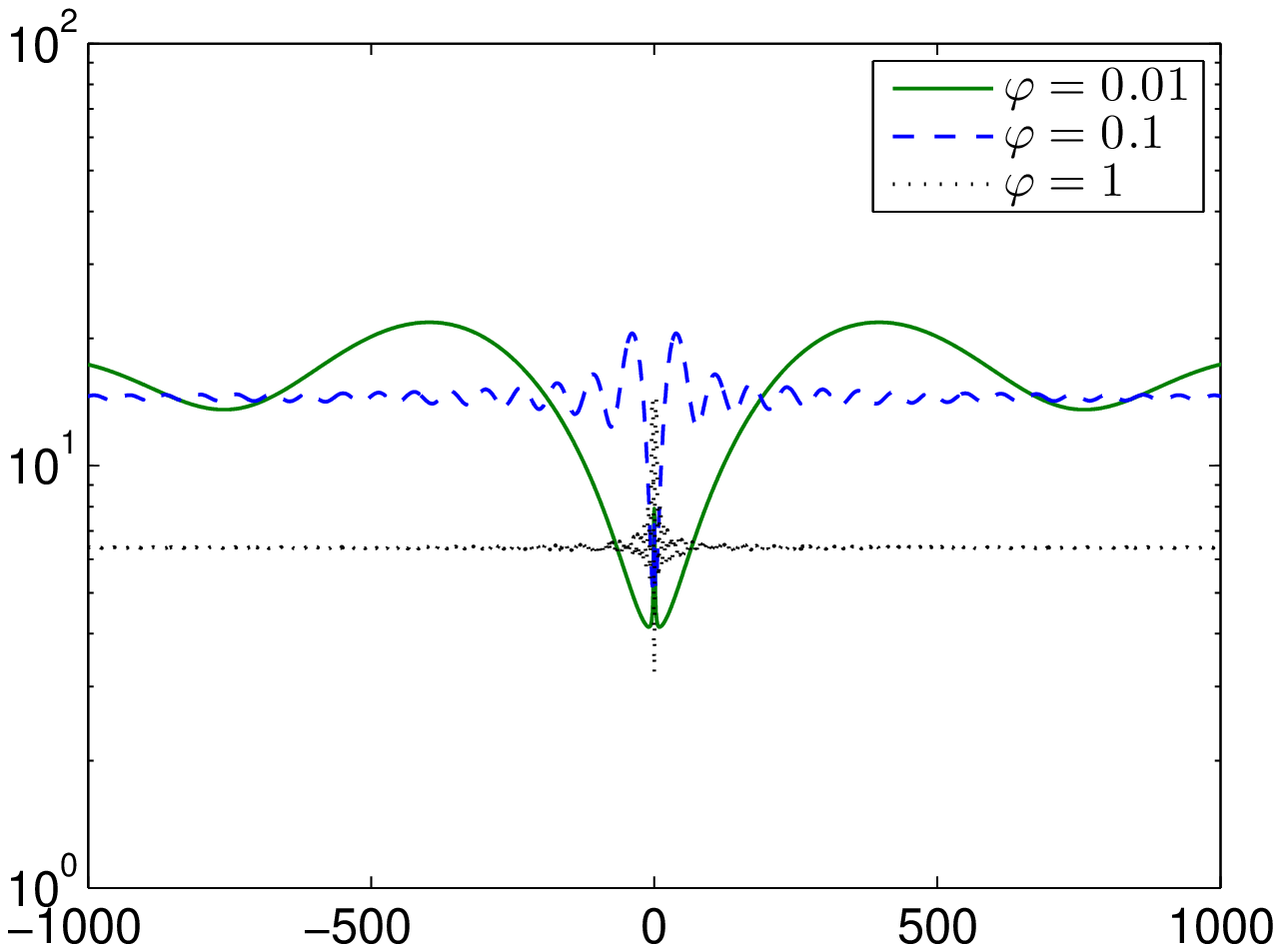}\includegraphics [scale=0.50]{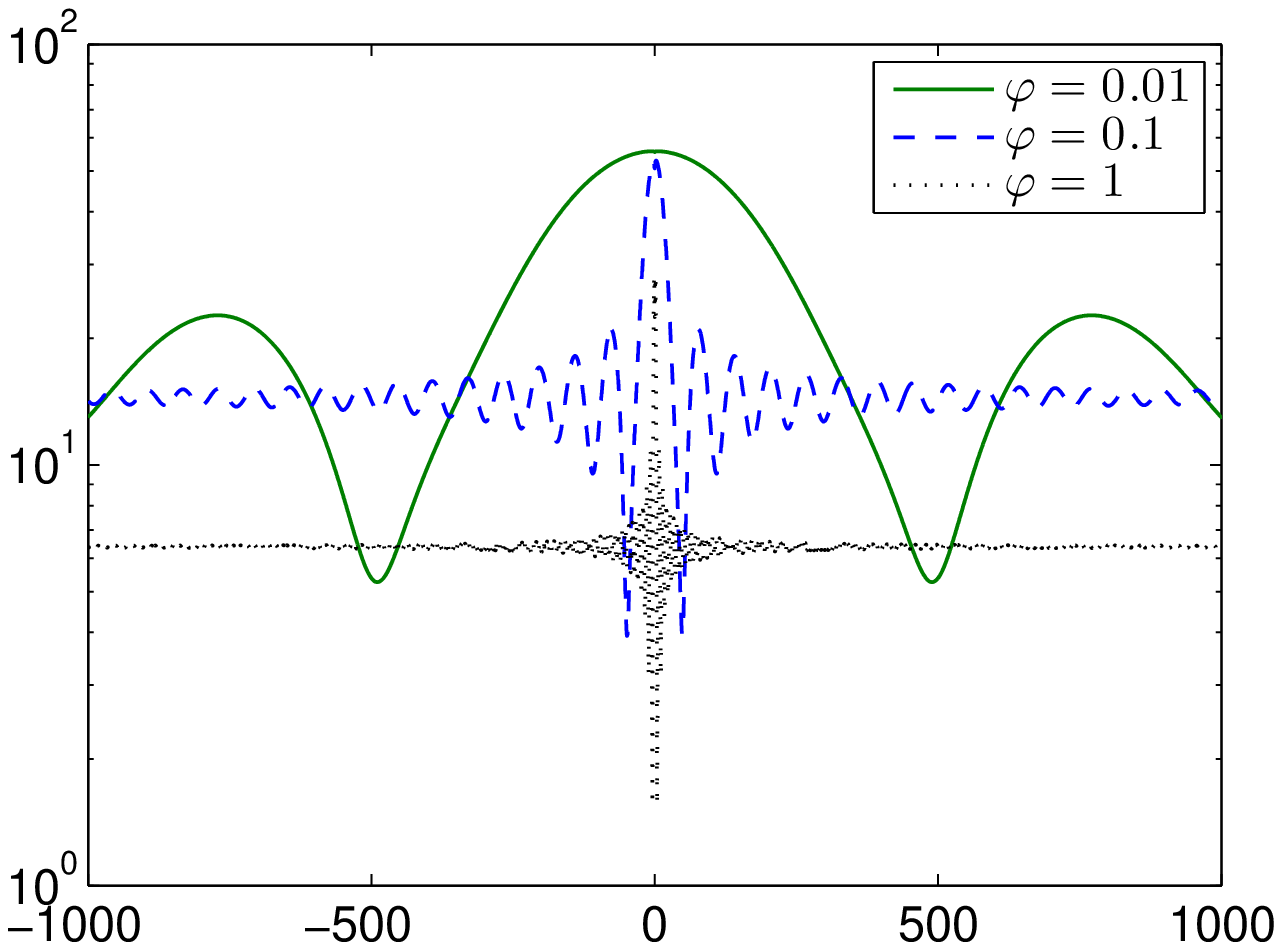}
       \put(-100,-5){\small$x$}
       \put(-390,125){$\varphi^{-2}m_{11}$}
       \put(-180,125){$\varphi^{-2}m_{22}$}
       \put(-320,-5){\small$x$}

  \end{center}
    \caption{\small Diagonal elements $m_{jj}(x)$ of the reminder matrix $M_{1,\varphi}(x)$ to be factorised on the next step for various values of the parameter $\varphi$. The elements are normalised to the value of the
    parameter $\varphi^2$. }
\label{f6}
\end{figure}

\begin{figure}[h!]
  \begin{center}
    \includegraphics [scale=0.50]{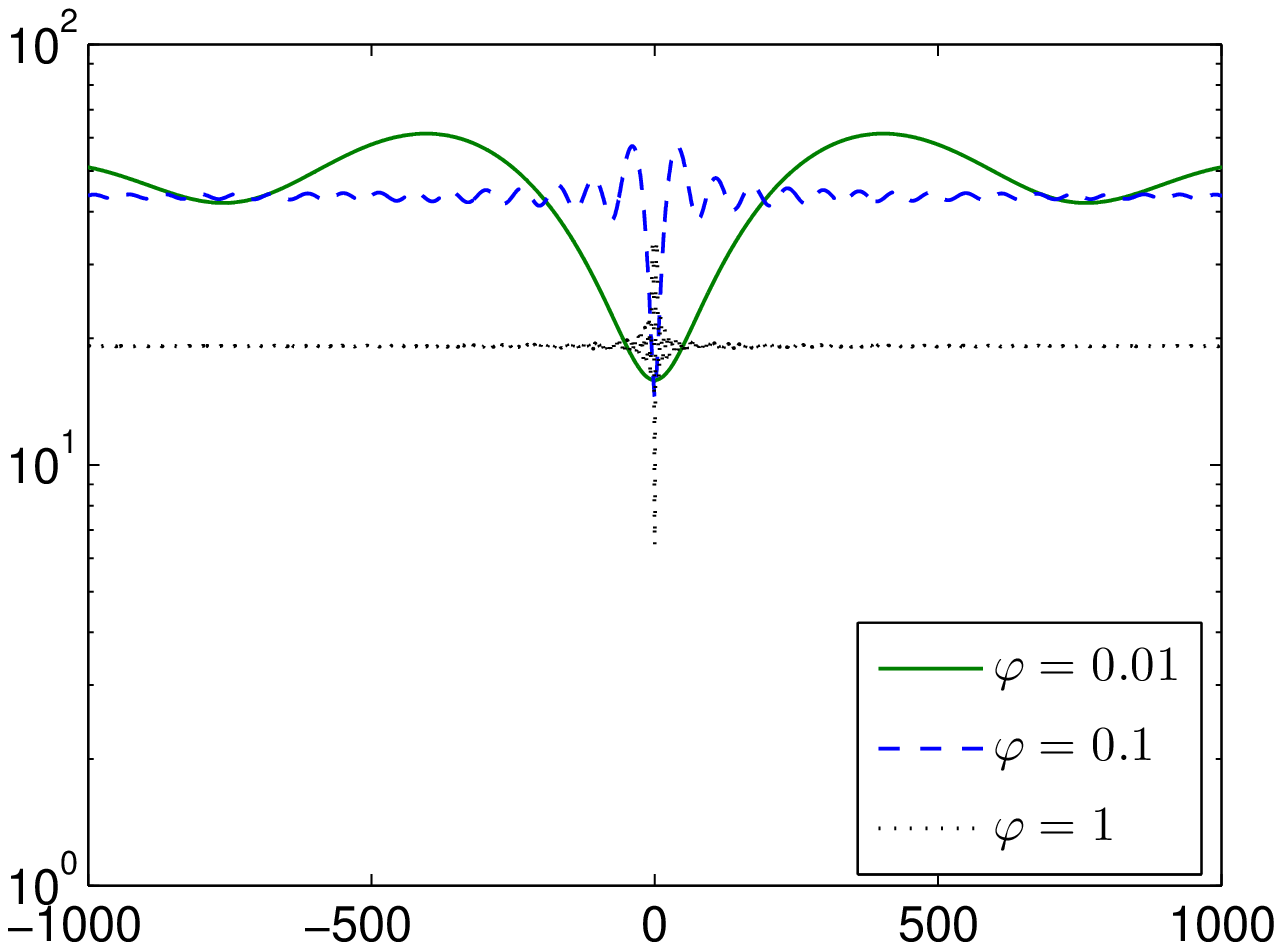}\includegraphics [scale=0.50]{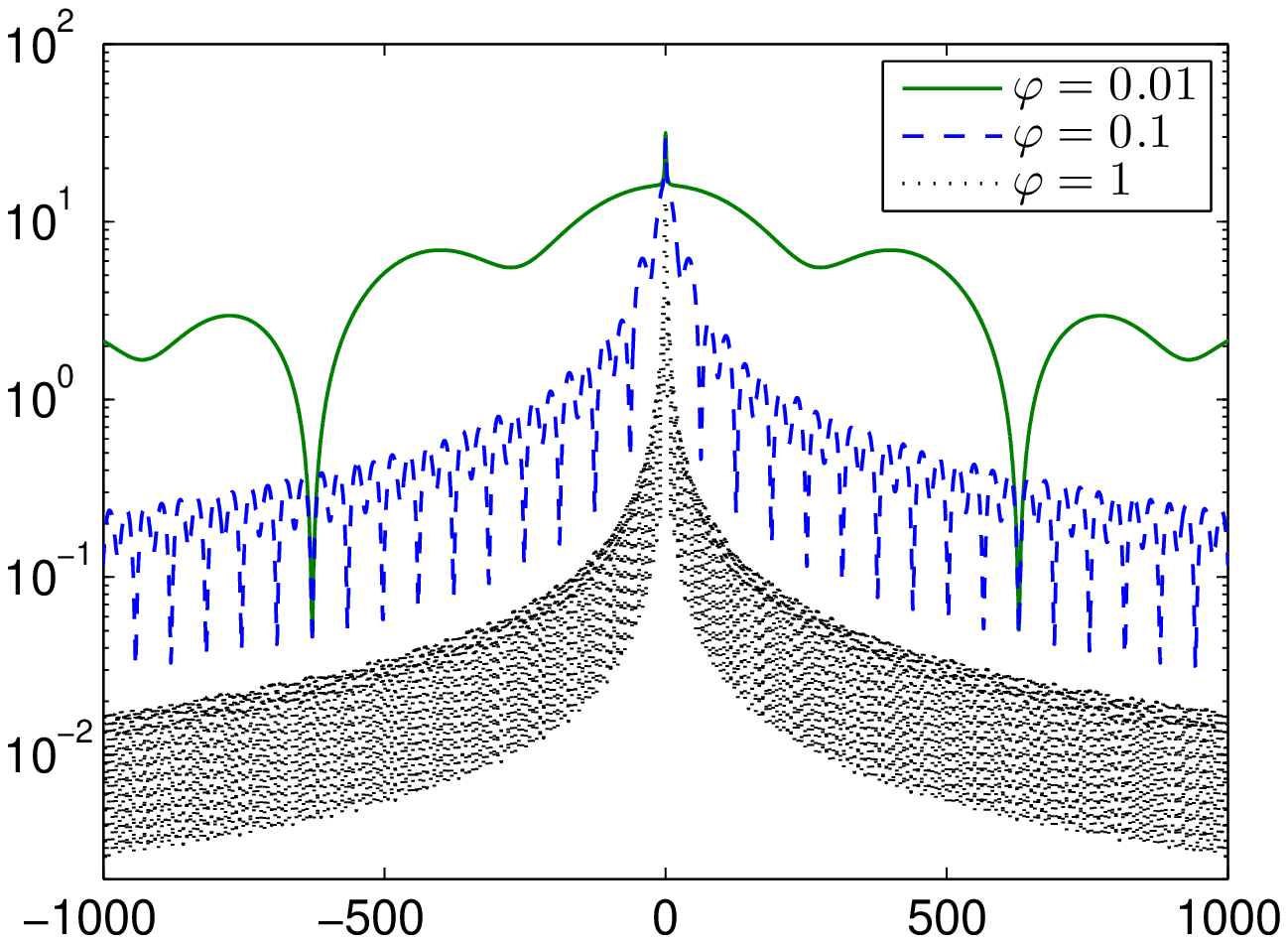}
       \put(-100,-5){\small$x$}
       \put(-390,30){$\varphi^{-2}m_{12}$}
       \put(-180,125){$\varphi^{-2}m_{21}$}
       \put(-320,-5){\small$x$}

  \end{center}
    \caption{\small The other two elements $m_{ij}(x)$ of the reminder matrix $M_{1,\varphi}(x)$ for various values of the parameter $\varphi$.
    The elements are normalised to the value of the parameter $\varphi^2$. }
\label{f7}
\end{figure}

In Fig.~\ref{f6} and Fig.~\ref{f7} we present results for this particular asymptotic expansion similar to those in Fig.~\ref{f1} and Fig.~\ref{f2}.
Now, apart of the fact we did not effect the max norms of the first asymptotic terms, $N^\pm_{1,\varphi}(x)$, the max norms of the reminder $M_{1,\varphi}(x)$  increased even further what is the price for the chosen matrix structure. Thus the convergence may become slow and the convergence radius for the small parameter $\varphi$
decreases on an order at least. As one can expect the maximal component of the reminder is $m_{22}$.

\subsubsection{Estimate of the reminder for the last choice of constant $C_0=C_0^{(3)}$}

The last case we would like to discuss refers to the choice of the constant matrix $C_0$ given in (\ref{Matrix_C_12})$_3$ which immediately leads to (\ref{Matrix_M_12})$_3$).
This, in this particular case, the components of the reminder, $M_{1,\varphi}(x)$, vanish at infinity similarly as it was for the original matrix function $M_{0,\varphi}(x)$.
As the results, the second step of the factorization is repetition of the first one in terms of its practical implementation
(e.i. there is no need to deal with the singular integral with the density having equal finite values at infinity (of the class $H_\mu(\overline{\mathbb{R}})$).
Clearly, such choice of the constant $C_0$ looks the most advantageous one from all the four examples considered in this paper. However, it is unlikely that such case is
available for an arbitrary factorization. An interesting observation is that the first and the last examples are very similar in terms
of the choice of the components of the constant $C_0$ (compare (\ref{Matrix_C_0}) and (\ref{Matrix_C_12})$_3$) and provide the most close values of the reminder $M_{1,\varphi}(x)$
in terms of the max norm. From comparing those examples, one can conclude that there is no much difference in terms of the asymptotic procedure to follow one or the others.
Even if some decrease in the convergence rate may follow from a choice it may sense to proceed in this way if the desired property is guaranteed at the end of the procedure.
Apart of the fact we can  not prove this, judging from the bahaviour of the reminder to be factorized on the next stage, the standard choice (\ref{Matrix_C_0})
of the constant $C_0$ provide the faster converging rate. The comparison also suggests that the minimization of the reminder norm is probably the most efficient strategy for the factorization. However, this restricts user's change to deliver the factorization with some specific properties of the factors.

\begin{figure}[h!]
  \begin{center}
    \includegraphics [scale=0.45]{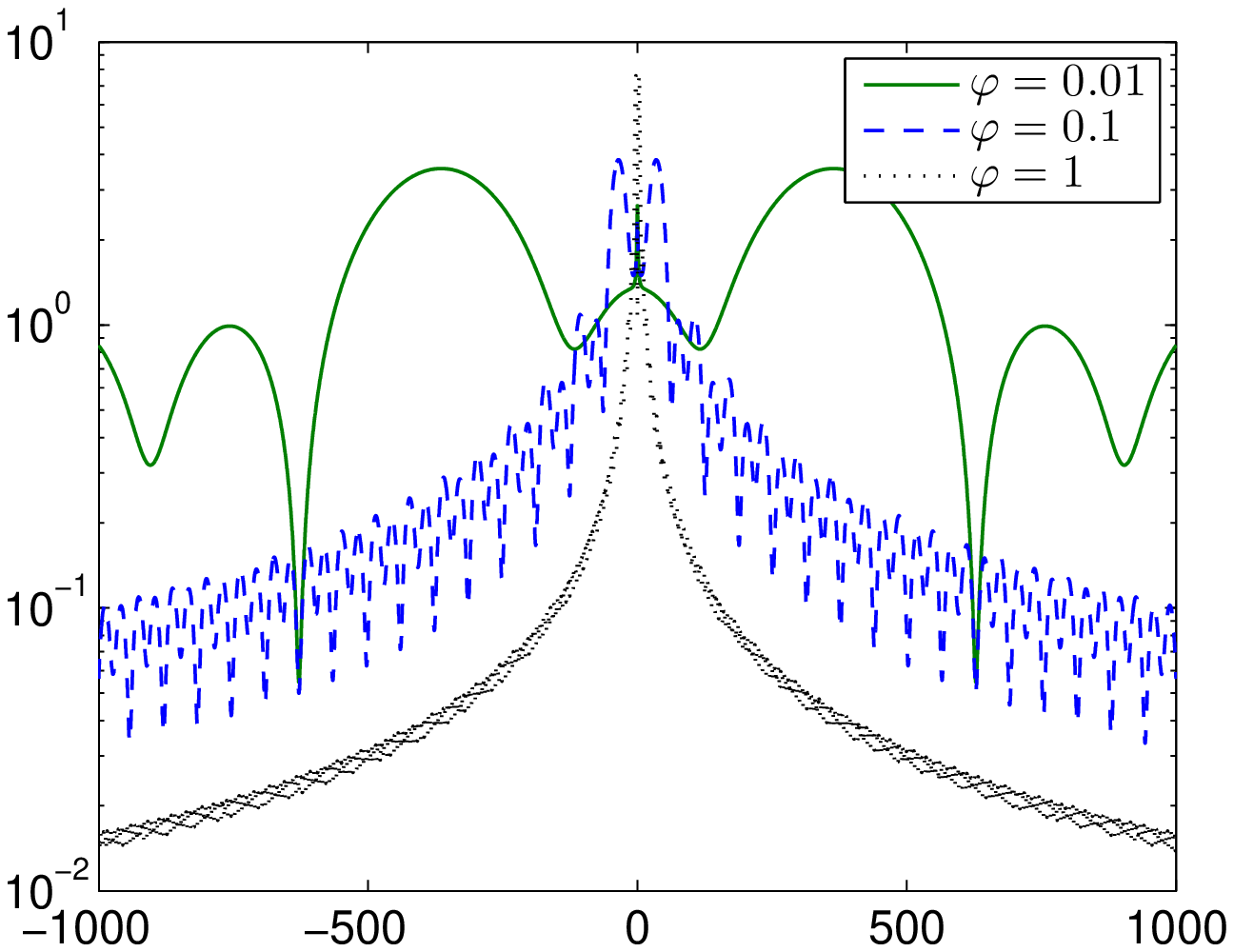}\includegraphics [scale=0.45]{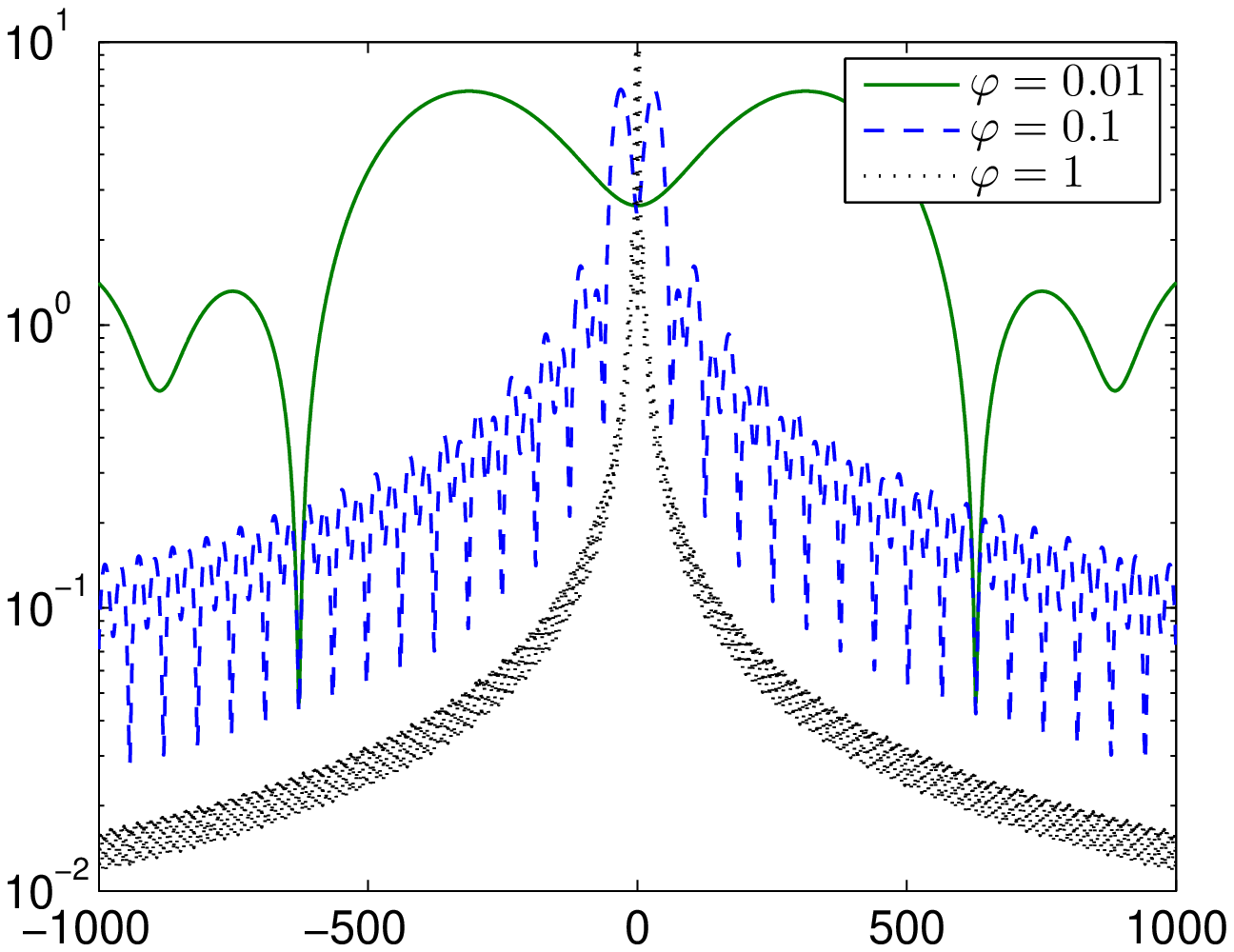}
       \put(-95,-5){\small$x$}
       \put(-300,20){$\varphi^{-2}m_{11}$}
       \put(-110,20){$\varphi^{-2}m_{22}$}
       \put(-285,-5){\small$x$}

  \end{center}
    \caption{\small Diagonal elements $m_{jj}(x)$ of the reminder matrix $M_{1,\varphi}(x)$ to be factorised on the next step for various values of the parameter $\varphi$. The elements are normalised to the value of the
    parameter $\varphi^2$. }
\label{f8}
\end{figure}

\begin{figure}[h!]
  \begin{center}
    \includegraphics [scale=0.45]{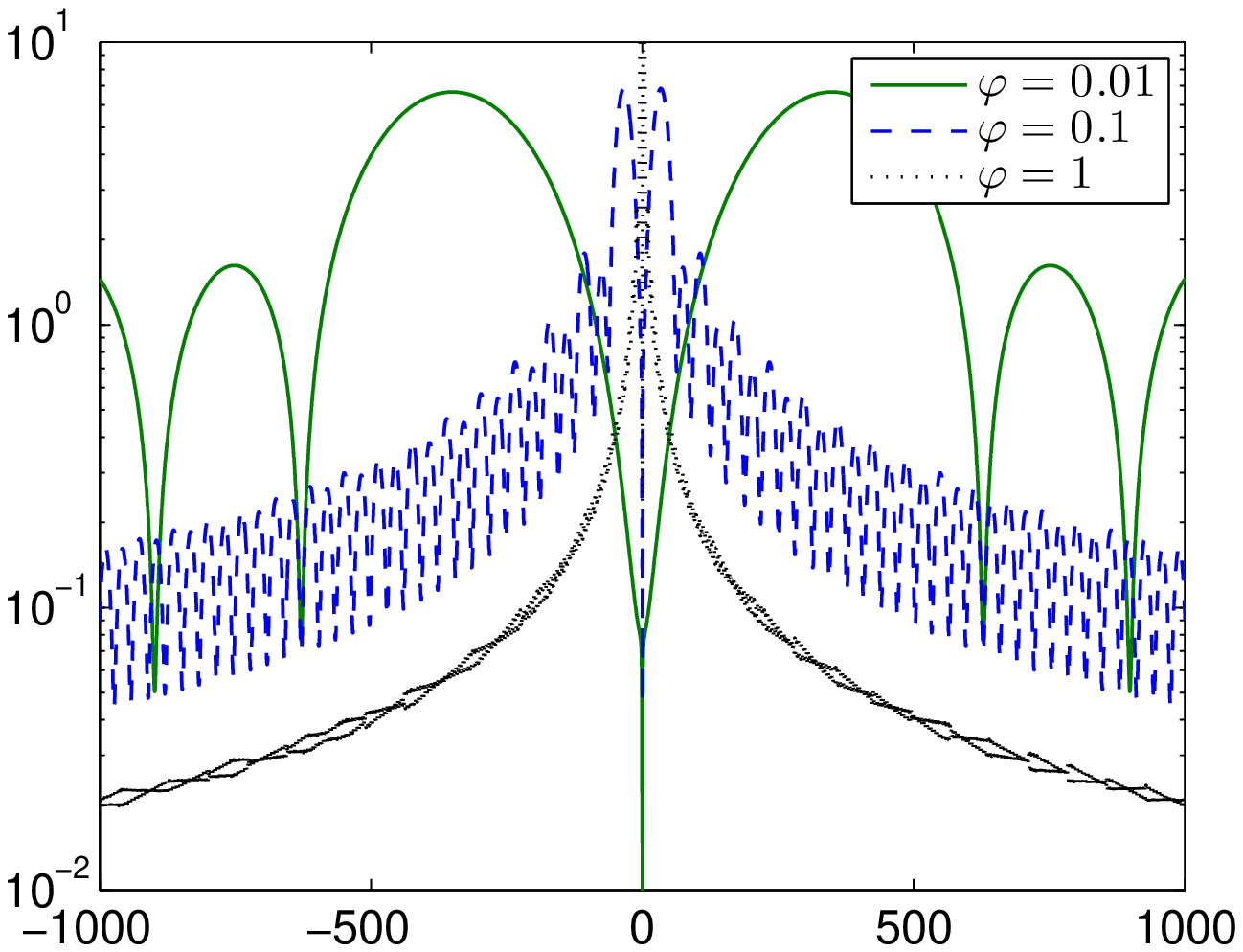}\includegraphics [scale=0.45]{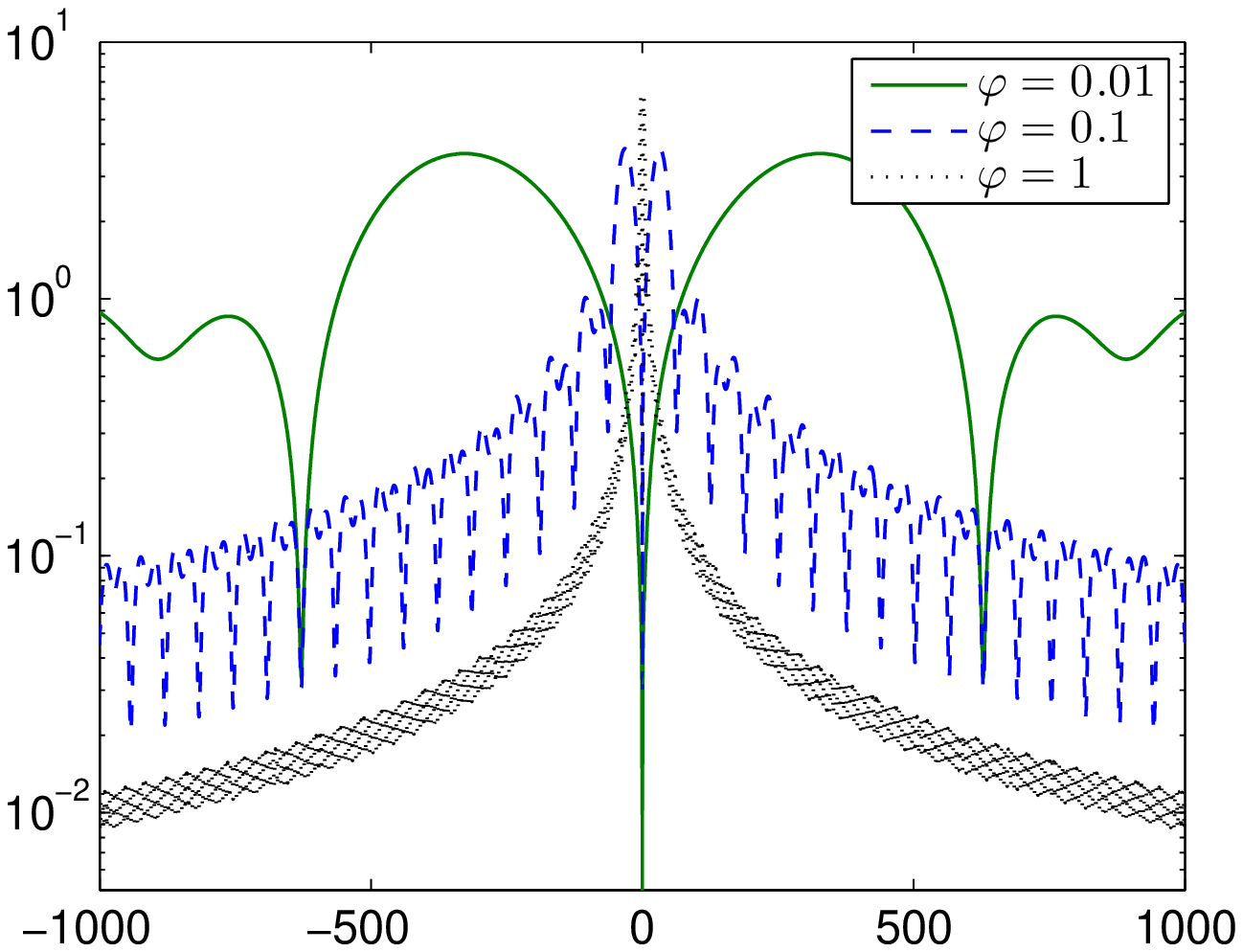}
       \put(-95,-5){\small$x$}
       \put(-300,20){$\varphi^{-2}m_{12}$}
       \put(-110,20){$\varphi^{-2}m_{21}$}
       \put(-285,-5){\small$x$}

  \end{center}
    \caption{\small The other two elements $m_{ij}(x)$ of the reminder matrix $M_{1,\varphi}(x)$ for various values of the parameter $\varphi$.
    The elements are normalised to the value of the parameter $\varphi^2$. }
\label{f9}
\end{figure}


\section{Conclusions}

We proposed the algorithm for effective asymptotic factorization of a class of nonrational matrix functions appeared in applications. We prove the convergence of the procedure thus, when the value of the small parameter lies in the convergence range, this factorization is explicitly represented in terms of series. The fact that the original unreturned matrix ($\varphi=0$) possesses the stable set of the partial indices is crucial for the method.  Although only right-hand factorization is considered here, the left-hand side analog follows trivially. The procedure does not lead to an unique form of the factorization factors and does not allow fixing factors limiting values at infinity by a given constant matrix. However, the limiting values of the factors for any factorization are represented by the nonsingular mutually invertible matrices (see (\ref{as_factors})$_2$). We show effectiveness of the process on a numerical example where even the first step provides a very good accuracy while the next steps of the procedure are easy repeatable. Finally, by choosing arbitrary constants
in a desirable way, one can direct the factorization process to preserve specific properties of the factors or to speed up the series convergence.


{\bf Acknowledgement.} {This research was supported by FP7-PEOPLE-2012-IAPP through the grant PIAP-GA-2012-284544-PARM2.}

\end{document}